\documentclass[10pt]{article}
\usepackage{amsmath}
\usepackage{latexsym}
\usepackage{amssymb}
\usepackage{amsfonts}
\usepackage{amsthm}
\title{ GRAVITATIONAL WAVES AND POMMARET BASES     }
\author{J.-F. POMMARET \\ CERMICS, Ecole des Ponts ParisTech, France \\ 
 jean-francois.pommaret@wanadoo.fr \\
ORCID: 0000-0003-0907-2601 }
\date{  }
\begin{document}\noindent
\maketitle

\noindent
{\bf ABSTRACT}  \\

\noindent
The first finite length differential sequence, now called {\it Janet sequence},  has been introduced 
by Janet in 1920. Thanks to the first book of Pommaret in 1978, this algorithmic approach has been extended by Gerdt, Blinkov, Zharkov, Seiler and others who introduced Janet and Pommaret bases in computer algebra. Between 1920 and 1990, new intrinsic tools have been developed by Spencer and successors in homological algebra, culminating with the definition of {\it extension differential modules} through the systematic use of {\it double differential duality}. If an operator ${\cal{D}}_1$ generates the compatibility conditions (CC) of an operator ${\cal{D}}$, then the {\it adjoint operator} $ad( {\cal{D}})$ may not generate the CC of $ad({\cal{D}}_1)$. Equivalently, an operator ${\cal{D}}$ with coefficients in a differential field $K$ can be parametrized by an operator ${\cal{D}}_{-1}$ iff the differential module $M$ defined by ${\cal{D}}$ is torsion-free, that is $t(M)={ext}^1_D(N,D) = 0$ when $N$ is the differential module defined by $ad( {\cal{D}})$ and $D$ is the ring of differential operators with coefficients in $K$. Also $R = hom_K(M,K)$ is a differential module for the Spencer operator $d:R \rightarrow T^* \otimes R$, first introduced by Macaulay in 1916 with {\it inverse systems}. When ${\cal{D}}$ is the self-adjoint Einstein operator, it is not evident that $t(M)\neq 0$ is generated by the components of the Weyl tensor having only to do with the pseudogroup of conformal transformations. Gravitational waves are thus not coherent with these results because the stress-functions parametrizing the Cauchy = ad (Killing) operator have nothing to do with the metric in general, {\it exactly like the Airy or Maxwell functions in elasticity}. Similarly, the Cauchy operator has nothing to do with any contraction of the Bianchi operator. The main difficulty is that $ad({\cal{D}})$ may not be  involutive {\it at all} when ${\cal{D}}$ is involutive, a well known fact in OD control theory leading to the Kalman test (Zbl 1079.93001). \\

\vspace{3cm}

\noindent
{\bf KEY WORDS}  \\

\noindent
Differential sequence; Spencer operator; Differential modules; ; Differential duality; \\
Extension modules; Control theory.  \\

\newpage

\noindent
{\bf 1) INTRODUCTION}. \\

Let us start this paper with a short personal but meaningful story that has oriented my research work during the last fifty years or so. In the fall of 1969 I decided to become a visiting student of D. C. Spencer in Princeton university, being attracted by learning his work at the source for future applications to physics. By chance, Spencer gave me his own key of the mathematical library, opened day and night and well furnished with french mathematical literature. However, if on one side I discovered that the intrinsic homological procedure developed by Spencer was exactly what I was dreaming about and decided to "bet" my life on it, on the other side it has been a very bad moment when I discovered that {\it Spencer and collaborators were totally unable to compute any explicit example}, the reason for which I had never found any one of them in their papers or books. The reader needs less than five minutes to discover that the introductory examples in the book ([17]) have not a single link with the core of this book. It is at this moment that I discovered, during a night in the library, the work of M. Janet written in $1920$ ([15]) that provided me the "Janet tabular" and the way to mix up a combinatoric approach with an intrinsic framework that led me to my first GB book appeared in $1978$ ([21]). This book has been translated into Russian by MIR in $1983$ with a successful distribution in what was called East of Europe because it was new and only costing $8$ roubles (about 5 US dollars at that time). This has been the origin of my first private contacts with V. Gerdt (Lectures in Moscow, Doubna and Iaroslav, 14-28/10/1995) before he introduced with Y. Blinkov what they called "Pommaret bases". The problem is that the computer algebra community did not understand that Spencer wanted to apply his methods for studying Lie pseudogroups, not at all for dealing with computers (See the Introduction of my first Kluwer book of 1994 ([24]) for the first computer study of an example provided by Janet). During many years, I tried vainly to convince Gerdt and the people of Aachen who were regularly inviting me that {\it the important side for applications is the intrinsic one}, even if {\it intrinsicness is always competing with complexity in computer algebra}. I gave it up after they supervised the thesis published in the reference [5] of ([46]) that must be compared to ([21]) with no need for comments. Indeed, one thing is to quote a book, another thing is to read it. In the meantime, I had the chance to meet Janet many times as he was living in Paris only a few blocks away from my parents and I can claim that his goal has always been to construct differential sequences along the footnote of his main $1920$ paper ([15]). The main purpose of this paper is to revisit the definition of Pommaret bases in the light of the differential sequence called "{\it Janet sequence}" and to explain why gravitational waves are not coherent with the results obtained. I do not know any other reference on the application of "{\it double differential duality} " to mathematical physics, the main difficulty being that the adjoint of an involutive operator may not be involutive at all for both OD and PD equations.  \\

With standard notations of differential geometry, let $(E, E', F, ...)$ be vector bundles over a manifold $X$ of dimension $n$ with sections $(\xi, \eta, ... )$, tangent bundle $T$ and cotangent bundle $T^*$. We shall denote by $J_q(E) $ the $q$-jet bundle of $E$ with sections ${\xi}_q$ transforming like the $q$-derivatives $j_q(\xi)$. If $\Phi: J_q(E) \rightarrow E'$ is a bundle morphism, we shall consider the system $R_q = ker(\Phi) \subset J_q(E)$ of order $q$ on $E$. The $r$-prolongation ${\rho}_r(R_q)= J_r(R_q) \cap J_{q+r}(E) \subset  J_r(J_q(E))$ obtained by differentiating formally $r$ times the given  ordinary (OD) or partial (PD) defining equations of $R_q$ will be the kernel of the composite morphism 
${\rho}_r(\Phi): J_{q++r}(E) \rightarrow J_r(J_q(E)) \stackrel{J_r(\Phi)}{\rightarrow} J_r(E')$. We shall also set $R'_r = im ({\rho}_r(\Phi)) \subset J_r(E')$. The symbol $g_{q+r} = R_{q+r} \cap S_{q+r}T^* \otimes E  \subset J_{q+r}(E) $ of $R_{q+r}$ is the $r$-prolongation of the symbol $g_q$ of $R_q$ and the kernel of the composite induced morphism ${\sigma}_r(\Phi): S_{q+r}T^* \otimes E \rightarrow S_rT^* \otimes E'$ obtained by restriction. We shall define the Spencer operator $d: R_{q+1} \rightarrow T^* \otimes R_q: {\xi}_{q+1} \rightarrow j_1 ({\xi}_q) - {\xi}_{q+1}$ by using the fact that $R_{q+1} = J_1(R_q) \cap J_{q+1}(E)$ and that $J_1(R_q)$ is an affine vector bundle over $R_q$ modelled on $T^* \otimes R_q$. We have the restriction $d:R'_{r+1} \rightarrow T^* \otimes R'_r$ and thus $R'_{r+1} \subset {\rho}_1(R'_r)$ . We shall suppose that $\Phi$ is an epimorphism and introduce the vector bundle $F_0 = J_q(E)/ R_q$. The system $R_q$ is said to be {\it formally integrable} (FI) if $r+1$ prolongations do not bring new equations of order $q+r$ other than the ones obtained after only $r$ prolongations, for any $r\geq 0$, that is all the equations of order $q+r$ can be obtained by differentiating $r$ times {\it only} the given equations of order $q$ for any $r \geq 0$. $R_q$ is said to be {\it involutive}  if it is FI and $g_q$ is involutive, a purely algebraic property [21, 24]. In that case, the successive CC operators can only be {\it at most} ${\cal{D}}_1, ... , {\cal{D}}_n$ which are first order involutive operators. \\

When $R_q$ is not involutive, a standard {\it prolongation/projection} (PP) procedure allows in general to find integers $r,s$ such that the image $R^{(s)}_{q+r}$ of the projection at order $q+r$ of the prolongation ${\rho}_{r+s}(R_q) = J_{r+s}(R_q) \cap J_{q+r+s}(E)\subset J_{r+s}(J_q(E)) $ is involutive with ${\rho}_t(R^{(s)}_{q+r})=R^{(s)}_{q+r+t}, \forall t\geq 0$ but it may highly depend on the parameters ([14, 21, 24]).  \\

The next problem is to define the CC operator ${\cal{D}}_1: F_0 \rightarrow F_1: \eta \rightarrow \zeta$ in such a way that the CC of ${\cal{D}} \xi = \eta$ is of the form 
${\cal{D}}_1 \eta = 0$. As shown in many books [21-24, 27, 37, 38] and papers [31, 39], such a problem may be quite difficult because the order of the generating CC may be quite high. Proceeding in this way, we may construct the CC ${\cal{D}}_2: F_1 \rightarrow F_2$ of ${\cal{D}}_1$ and so on. The difficulty, shown on the motivating examples, is that " {\it jumps} " in the successive orders may appear, even on elementary examples. Now, if the map $\Phi$ depends on constant (or variable) parameters $(a, b, c, ... )$, then the study of the two previous problems becomes much harder because the ranks of the matrices ${\rho}_r(\Phi)$ and/or ${\sigma}_r(\Phi)$ may also highly depend on the parameters as we shall see. Such a question is particularly delicate in the study of the Kerr $(m, a)$, Schwarzschild $(m, 0)$ and Minkowski $(0, 0)$ parameters while computing the dimensions of the inclusions 
$R^{(3)}_1\subset R^{(2)}_1 \subset R^{(1)}_1 =R_1 \subset J_1(T)$ for the respective Killing operators as the numbers of generating second order CC and the numbers of generating third order CC may change drastically [45].   \\ 

\noindent
{\bf EXAMPLE 1.1}: Let $n=2, m= 1$ and introduce the trivial vector bundle $E$ with local coordinates $(x^1, x^2, \xi)$ for a section over the base manifold $X$ with local coordinates $(x^1,x^2)$. Let us consider the linear second order system $R_2 \subset J_2(E)$ defined by the two linearly independent equations $ d_{22}\xi= 0, \,\, d_{12}\xi + a d_1 \xi = 0$ where $a$ is an arbitrary constant parameter. Using crossed derivatives, we get the second order system $R^{(1)}_2 \subset R_2$ defined by the PD equations 
$d_{22} \xi = 0, d_{12} \xi + a d_1 \xi= 0, a^2 d_1\xi=0$ 
which is easily seen not to be involutive. Hence we have two possibilities:  \\
 $\bullet  \,\,\, a=0 $: We obtain the following second order homogeneous involutive system: \\
\[ R^{(1)}_2 = R_2 \subset J_2(E) \hspace{2cm} \left\{  \begin{array}{rcl}
 d_{22} \xi &=  & {\eta}^2 \\
d_{12} \xi  & = & {\eta}^1  
\end{array}
\right.  \fbox{ $ \begin{array}{ll}
1 & 2   \\
1 & \bullet 
\end{array} $ } \] 
with only one first order CC operator $d_2 {\eta}^1 - d_1 {\eta}^2= \zeta $,  leading to the Janet sequence:  \\
\[   0 \rightarrow  \Theta \rightarrow E \underset 2 {\stackrel{{\cal{D}}}{\longrightarrow }} F_0  \underset 1 {\stackrel{{\cal{D}}_1}{\longrightarrow }} F_1 \rightarrow 0 
 \hspace{2cm}   \xi  \stackrel{{\cal{D}}}{\longrightarrow } \eta  \stackrel{{\cal{D}}_1}{\longrightarrow }  \zeta  \rightarrow 0 \]
Multiplying by $\lambda$ and integrating by parts, we get the operator $ad({\cal{D}}_1)$ which is described by $ - d_2 \lambda = {\mu}^1, -d_1 \lambda = {\mu}^2$.
Multiplying the first equation by ${\mu}^1$, the second by ${\mu}^2$, summing and integrating by parts, we notice that $ad({\cal{D}})$, described by $d_{12}{\mu}^1 + d_{22} {\mu}^2 =\nu$, is of order $2$ and does not therefore generates the CC of $ad({\cal{D}}_1)$ which is of order $1$, namely $d_1 {\mu}^1 + d_2{\mu}^2 = {\nu}'$ as below: \\
\[        \begin{array}{rcccl}
 \nu & \stackrel{ad({\cal{D}})}{\longleftarrow} &  \mu  & \stackrel{ad({\cal{D}}_1)}{\longleftarrow} &  \lambda   \\
  &  \swarrow &  &  &  \\
  {\nu}' &  &  &  &  
 \end{array}  \]
 and $ {\nu}'$ is a torsion element of the differential module defined by $ad({\cal{D}})$ because $d_2 {\nu}' = \nu$. As we shall see in the third section, if $M_1$ is the differential module defined by ${\cal{D}}_1$, then we know that ${ext}^1(M_1) \neq 0$ when $a=0$.  \\
 
\noindent
$\bullet  \,\,\, a \neq 0$: We obtain the second order system $R^{(1)}_2 $ defined by $ d_{22} \xi = 0, d_{12} \xi = 0, d_1\xi = 0 $ with a strict inclusion $R^{(1)}_2 \subset R_2$ because $3<4$. We may define $\eta = d_1{\eta}^2 - d_2 {\eta}^1+ a {\eta}^1$ and obtain the involutive and finite type system in $\delta$-regular coordinates:  \\
\noindent
\[ R^{(2)}_2 \subset J_2(E) \hspace{2cm} \left\{  \begin{array}{rcl}
 d_{22} \xi &=  & {\eta}^2  \\
 d_{12} \xi & =  & {\eta}^1 - \frac{1}{a} \eta   \\
 d_{11} \xi & =  & \frac{1}{a^2} d_1 \eta \\
 d_1  \xi    & = & \frac{1}{a^2} \eta
\end{array}
\right.  \fbox{ $ \begin{array}{ll}
1 & 2   \\
1 & \bullet \\  
1 & \bullet  \\
\bullet & \bullet 
\end{array} $ } \]
Counting the dimensions, we have the following strict inclusions by comparing the dimensions:  \\
\[      R^{(2)}_2 \subset R^{1)}_2 \subset R_2 \subset J_2(E) , \hspace{3cm}  2 <  3   <  4  <  6  \]

The symbol $g_{2+r} $ is involutive with $dim(g_{2+r})=1, \forall r\geq 0$ and we have 
$dim(R_{2+r})  = 4, \forall r\geq 0$. Indeed, using jet notation, the $4$ parametric jets of $R_2$ are $(\xi, {\xi}_1, {\xi}_2, {\xi}_{11})$. The $4$ parametric jets of $R_3$ are now $(\xi, {\xi}_2, {\xi}_{11}, {\xi}_{111})$ and so on. Accordingly, the dimension of $g_{2+r}$ is $1$ because the only parametric jet is ${\xi}_{1....1}$. We have the short exact sequence 
$0 \rightarrow g_{r+2} \rightarrow R_{r+2} \rightarrow R^{(1)}_{r+1} \rightarrow 0 $ and the symbol of $R_2$ is involutive. It follows from a delicate but crucial theorem (See [21], Theorem 2.4.5 p 70 and proposition 2.5.1 p 76 with $q=2$) that ${\rho}_r(R^{(1)}_2) = R^{(1)}_{2+r} \,\,, \forall r\geq 0$ with $dim(R^{(1)}_{r+1})= 3$, a result leading to $dim(R_{r+2})= 3+1=4$ by counting the dimensions. As $R^{(1)}_2$ does not depend any longer on the parameter, the general solution is easily seen to be of the form $\xi = cx^2 + d$ and is thus only depending on two arbitrary constants, contrary to what could be imagined from this result but in a coherent way with the fact that $dim(R^{(2)}_{2+r}) = 2, \forall r\geq 0$.    \\

After differentiating twice, we could be waiting for CC of order $3$. However, we obtain the $4$ CC:  \\
\[ d_2 {\eta}^1 - \frac{1}{a} d_2 \eta - d_1 {\eta}^2 = 0, \frac{1}{a^2} d_{12} \eta - d_1 {\eta}^1 + \frac{1}{a} d_1\eta = 0, \frac{1}{a^2} d_2 \eta - {\eta}^1 + \frac{1}{a} \eta =0, 
\frac{1}{a^2} (d_1 \eta - d_1 \eta) = 0  \]
The last CC that we shall call "{\it identity to zero} " must not be taking into account. The second CC is just the derivative with respect to $x^1 $ of the third CC which amounts to 
\[ (d_{12} {\eta}^2 - d_{22}{\eta}^1 + a d_2 {\eta}^1) - a^2 {\eta}^1 + a (d_1 {\eta}^2 - d_2 {\eta}^1 + a {\eta}^1)=0 \Leftrightarrow  d_{12} {\eta}^2 - d_{22} {\eta}^1 + a d_1 {\eta}^2=0 \]
which is a second order CC amounting to the first. Hence we get the only generating CC operator ${\cal{D}}_1: ({\eta}^1, {\eta}^2) \rightarrow d_{12} {\eta}^2 - d_{22} {\eta}^1 + a d_1 {\eta}^2=\zeta $ which is thus formally surjective.  \\

For helping the reader, we recall that basic elementary combinatorics arguments are giving $dim(S_qT^*) = q+1$ while $dim(J_q(E)) = (q+1)(q+2)/2$ because $n=2$ and $m = dim(E)=1$. 
We obtain successively till we stop:
\[ R'_0=F_0  \Rightarrow R'_1=J_1(F_0) \Rightarrow  {\rho}_1(R'_1) = J_2(F_0) \Rightarrow R'_2 \subset {\rho}_1(R'_1) \Rightarrow  R'_3={\rho}_1(R'_2) \]
Hence, the number of generating CC of order $1$ is zero and the number of generating CC of strict order $2$ is $dim({\rho}_1(R'_1)) - dim(R'_2) = 12 - (15-4)=12-11=1$ in a 
coherent way.    \\
Setting $F_1 = Q_2$ with $dim(Q_2) = 1$, we obtain the commutative diagram:  \\

\[ \fbox{ $  \begin{array}{rcccccccl}
    &  0  &  &  0  &  &  0  &  &  0  &   \\
     &  \downarrow &  & \downarrow &  &  \downarrow  &  & \downarrow   &     \\
0 \rightarrow  &  g_5  & \rightarrow & S_5T^* \otimes E & \rightarrow &  \fbox { $ S_3 T^*\otimes F_0 $} & \rightarrow &  T^* \otimes F_1 & \rightarrow 0  \\
&  \downarrow &  & \downarrow &  &  \downarrow  &  &  \downarrow    &     \\
0 \rightarrow &  R_5 & \rightarrow  & J_5(E) & \rightarrow &  J_3(F_0) &  \rightarrow &  J_1(F_1) & \rightarrow 0   \\
    &  \downarrow &  & \downarrow &  &  \downarrow  &  & \downarrow   &     \\
0 \rightarrow &  R_4 & \rightarrow  & J_4(E) & \rightarrow &  J_2(F_0) &  \rightarrow &  F_1 &  \rightarrow 0             \\ 
   &  &  & \downarrow &  &  \downarrow  &  &  \downarrow   & \\
   &  &  &  0  && 0  && 0 &
\end{array} $ }  \]
with dimensions:
\[   \begin{array}{rcccccccl}
    &  0  &  &  0  &  &  0  &  &  0  &   \\
     &  \downarrow &  & \downarrow &  &  \downarrow  &  & \downarrow   &     \\
0 \rightarrow  &  1 & \rightarrow &  6  & \rightarrow &  \fbox { $  8 $} & \rightarrow &  2 & \rightarrow 0  \\
&  \downarrow &  & \downarrow &  &  \downarrow  &  &  \downarrow    &     \\
0 \rightarrow &  4 & \rightarrow  & 21 & \rightarrow & 20 &  \rightarrow &  3 & \rightarrow 0   \\
    &  \downarrow &  & \downarrow &  &  \downarrow  &  & \downarrow   &     \\
0 \rightarrow &  4 & \rightarrow  & 15 & \rightarrow &  12 &  \rightarrow &  1  &  \rightarrow 0             \\ 
   &  &  & \downarrow &  &  \downarrow  &  &  \downarrow   & \\
   &  &  &  0  && 0  && 0 &
\end{array}   \]
The upper symbol sequence is not exact at $S_3 T^* \otimes F_0$ even though the two other sequences are exact on the jet level. As a byproduct we have the exact sequences 
$\forall r\geq 0$:  \\
\[   0 \rightarrow  R_{r + 4} \rightarrow J_{r + 4}(E) \rightarrow  J_{r + 2}(F_0) \rightarrow  J_r(F_1) \rightarrow 0   \] 
Such a result can be checked directly through the identity:   \\
 \[   4 - (r+5)(r+6)/2 + 2(r+3)(r+4)/2 - (r+1)(r+2)/2 = 0   \]
We obtain therefore the formally exact sequence we were looking for, namely:  \\
\[  0 \rightarrow \Theta \rightarrow E \underset 2{\stackrel{{\cal{D}}}{\longrightarrow}}  F_0  \underset 2{\stackrel{{\cal{D}}_1}{\longrightarrow}}  F_1 \rightarrow 0   \]
{\it The surprising fact is that, in this case}, $ad({\cal{D}}) $ {\it generates the CC of} $ad({\cal{D}}_1) $. Indeed, multiplying by the Lagrange multiplier test function $\lambda$ and integrating by parts, we obtain the second order operator $ \lambda \rightarrow (- d_{22} \lambda = {\mu}^1, d_{12} \lambda - a d_1 \lambda = {\mu}^2 ) $  and thus 
$ - a^2 d_1 \lambda = d_1 {\mu}^1 + d_2 {\mu}^2  + a {\mu}^2 $. Substituting, we finally get the only second order CC operator $d_{12} {\mu}^1 + d_{22}{\mu}^2 - a d_1 {\mu}^1=0$. As we shall see in the third section, we have {\it now} ${ext}^1(M_1) = 0 $ when $a \neq 0$ and the adjoint sequences:
\[  \xi  \stackrel{{\cal{D}}}{\longrightarrow } \eta  \stackrel{{\cal{D}}_1}{\longrightarrow }  \zeta  \rightarrow 0 \hspace{2cm}
  0 \longleftarrow    \nu  \stackrel{ad({\cal{D}})}{\longleftarrow}   \mu   \stackrel{ad({\cal{D}}_1)}{\longleftarrow}   \lambda    \]                                                                                                

In the differential module framework over the commutative ring $ D= K[d_1,d_2]$ of differential operators with coefficients in the trivially differential field $K= \mathbb{Q}(a)$, we have the free resolution:
\[         0 \rightarrow  D  \underset 2{\stackrel{{\cal{D}}_1}{\longrightarrow}} D^2   \underset 2{\stackrel{{\cal{D}}}{\longrightarrow}} D \rightarrow M \rightarrow 0 \]
of the differential module $ M $ with Euler-Poincar\'{e} characteristic $rk_D(M) = 1 - 2 + 1 = 0$. We recall that $R= R_{\infty} = hom_K(M,K)$ is a differential module for the Spencer operator $d: R \rightarrow T^* \otimes R : R_{q+1}  \rightarrow T^* \otimes R_q$ (See [27, 38] for more details).   \\

Having in mind the pabove example, we have revisited these works by using new homological techniques and it is a matter of {\it fact} that they do not agree with the previous ones for the third order CC. In order to escape from such an unpleasant situation, we have written ([45]) in such a way that we are only using elementary combinatorics and diagram chasing. However, an equally important second purpose is to question the proper target of the quoted problem. Indeed, important concepts such as {\it differential extension modules} have been introduced in {\it differential homological algebra} and are known, thanks to a quite difficult theorem [20, 27, 52], to be the only intrinsic results that could be obtained independently of the differential sequence that could be used, provided that $\Theta= \ker({\cal{D}})$ is the same, that is even if one is using another system on $E$ with the same solutions. Equivalently, this amounts to say, in a few words but a more advanced language, if {\it we are keeping the same differential module} $ M $ {\it but changing its presentation}. \\

Of course, {\it in general} as we just saw, the extension modules may highly depend on the parameters. However, as we just saw, there are even simple academic systems depending on parameters but such that {\it a convenient equivalent system, say involutive with the same solutions, may no longer depend on the parameters} and the extension modules do not depend on the parameters because it is known that they do not depend on the differential sequence used for their definition.  This will be {\it exactly} the situation met in the study of the Kerr $(m, a)$, Schwarzschild $(m, 0)$ and Minkowski $(0, 0)$ parameters while studying the respective Killing operators ([45]). \\

Also, in a totally independent way which is still not acknowledged, E. Vessiot has shown that certain operators may depend on geometric objects satisfying non-linear {\it structure equations} that are depending on certain {\it Vessiot structure constants} $c$. The simplest example is the condition of constant Riemannian curvature [21-24, 37] which is necessary in order that the Killing system becomes FI but the case of classical or unimodular contact structures is similar ([46]). We have proved that the extension modules only depend on these constants ([39, 40, 42]). \\   \\

\noindent
{\bf 2) DIFFERENTIAL SYSTEMS}  \\

If $X$ is a manifold of dimension $n$ with local coordinates $(x)=(x^1, ... ,x^n)$, we denote as usual by $T=T(X)$ the {\it tangent bundle} of $X$, by $T^*=T^*(X)$ the {\it cotangent bundle}, by ${\wedge}^rT^*$ the {\it bundle of r-forms} and by $S_qT^*$ the {\it bundle of q-symmetric tensors}. More generally, let $E$ be a {\it vector bundle} over $X$ with local coordinates $(x^i,y^k)$ for $i=1,...,n$ and $k=1,...,m$ simply denoted by $(x,y)$, {\it projection} $\pi:E\rightarrow X:(x,y)\rightarrow (x)$ and changes of local coordinate $\bar{x}=\varphi (x), \bar{y}=A(x)y$. We shall denote by $E^*$ the vector bundle obtained by inverting the matrix $A$ of the changes of coordinates, exactly like $T^*$ is obtained from $T$. We denote by $f:X\rightarrow E: (x)\rightarrow (x,y=f(x))$ a global {\it section} of $E$, that is a map such that $\pi\circ f=id_X$ but local sections over an open set $U\subset X$ may also be considered when needed. Under a change of coordinates, a section transforms like $\bar{f}(\varphi(x))=A(x)f(x)$ and the changes of the derivatives can also be obtained with more work. We shall denote by $J_q(E)$ the {\it q-jet bundle} of $E$ with local coordinates $(x^i, y^k, y^k_i, y^k_{ij},...)=(x,y_q)$ called {\it jet coordinates} and sections $f_q:(x)\rightarrow (x,f^k(x), f^k_i(x), f^k_{ij}(x), ...)=(x,f_q(x))$ transforming like the sections $j_q(f):(x) \rightarrow (x,f^k(x), {\partial}_if^k(x), {\partial}_{ij}f^k(x), ...)=(x,j_q(f)(x))$ where both $f_q$ and $j_q(f)$ are over the section $f$ of $E$. For any $q\geq 0$, $J_q(E)$ is a vector bundle over $X$ with projection ${\pi}_q$ while $J_{q+r}(E)$ is a vector bundle over $J_q(E)$ with projection ${\pi}^{q+r}_q, \forall r\geq 0$.\\

Let $\mu=({\mu}_1,...,{\mu}_n)$ be a multi-index with {\it length} ${\mid}\mu{\mid}={\mu}_1+...+{\mu}_n$, {\it class} $i$ if ${\mu}_1=...={\mu}_{i-1}=0,{\mu}_i\neq 0$ and $\mu +1_i=({\mu}_1,...,{\mu}_{i-1},{\mu}_i +1, {\mu}_{i+1},...,{\mu}_n)$. We set $y_q=\{y^k_{\mu}{\mid} 1\leq k\leq m, 0\leq {\mid}\mu{\mid}\leq q\}$ with $y^k_{\mu}=y^k$ when ${\mid}\mu{\mid}=0$. If $E$ is a vector bundle over $X$ and $J_q(E)$ is the $q$-{\it jet bundle} of $E$, then both sections $f_q\in J_q(E)$ and $j_q(f)\in J_q(E)$ are over the section $f\in E$. There is a natural way to distinguish them by introducing the {\it Spencer}  operator $d:J_{q+1}(E)\rightarrow T^*\otimes J_q(E)$ with components $(df_{q+1})^k_{\mu,i}(x)={\partial}_if^k_{\mu}(x)-f^k_{\mu+1_i}(x)$. The kernel of $d$ consists of sections such that $f_{q+1}=j_1(f_q)=j_2(f_{q-1})=...=j_{q+1}(f)$. Finally, if $R_q\subset J_q(E)$ is a {\it system} of order $q$ on $E$ locally defined by linear equations ${\Phi}^{\tau}(x,y_q)\equiv a^{\tau\mu}_k(x)y^k_{\mu}=0$ and local coordinates $(x,z)$ for the parametric jets up to order $q$, the $r$-{\it prolongation} $R_{q+r}={\rho}_r(R_q)=J_r(R_q)\cap J_{q+r}(E)\subset J_r(J_q(E))$ is locally defined when $r=1$ by the linear equations ${\Phi}^{\tau}(x,y_q)=0, d_i{\Phi}^{\tau}(x,y_{q+1})\equiv a^{\tau\mu}_k(x)y^k_{\mu+1_i}+{\partial}_ia^{\tau\mu}_k(x)y^k_{\mu}=0$ and has {\it symbol} $g_{q+r}=R_{q+r}\cap S_{q+r}T^*\otimes E\subset J_{q+r}(E)$ if one looks at the {\it top order terms}. If $f_{q+1}\in R_{q+1}$ is over $f_q\in R_q$, differentiating the identity $a^{\tau\mu}_k(x)f^k_{\mu}(x)\equiv 0$ with respect to $x^i$ and substracting the identity $a^{\tau\mu}_k(x)f^k_{\mu+1_i}(x)+{\partial}_ia^{\tau\mu}_k(x)f^k_{\mu}(x)\equiv 0$, we obtain the identity $a^{\tau\mu}_k(x)({\partial}_if^k_{\mu}(x)-f^k_{\mu+1_i}(x))\equiv 0$ and thus the restriction $d:R_{q+1}\rightarrow T^*\otimes R_q$. More generally, we have the restriction:   \\
\[   d: {\wedge}^sT^* \otimes R_{q+1} \rightarrow {\wedge}^{s+1}T^* \otimes R_q: (f^k_{\mu,I}(x)dx^I) \rightarrow (({\partial}_if^k_{\mu,I}(x) - f^k_{\mu + 1_i,I}(x))dx^i \wedge dx^I)\] 
with standard multi-index notation for exterior forms and one can easily check that $d\circ d=0$. The restriction of $-d$ to the symbol is called the {\it Spencer} map $\delta$ in the sequences:  \\
\[ {\wedge}^{s-1}T^*\otimes g_{q+1} \stackrel{\delta}{\longrightarrow} {\wedge}^s T^* \otimes g_q 
\stackrel{\delta}{\longrightarrow} {\wedge}^{s+1} T^* \otimes g_{q-1}  \]
because $\delta \circ \delta=0$ similarly, leading to the purely algebraic {\it $\delta$-cohomology} $H^s_{q+r}(g_q)$ at $ {\wedge}^sT^* \otimes g_q$ [14, 21, 24, 27, 38, 54]. \\

\noindent
{\bf DEFINITION 2.1}: If $R_q\subset J_q(E)$ is a system of order $q$ on $E$, then $R_{q+r}={\rho}_r(R_q)=J_r(R_q) \cap J_{q+r}(E)\subset J_r(J_q(E))$ is called the {\it r-prolongation} of 
$R_q$. In actual practice, if the system is defined by PDE ${\Phi}^{\tau} \equiv a^{\tau \mu}_k (x) y^k_{\mu}=0$ the first prolongation is defined by adding the PDE $d_i{\Phi}^{\tau}\equiv a^{\tau \mu}_k(x) y^{\mu + 1_i} + {\partial}_i a^{\tau mu}_k(x)  y^k_{\mu} = 0$. Accordingly, $f_q \in R_q \Leftrightarrow a^{\tau \mu}_k(x) f^k_{\mu}(x) = 0$ and $f_{q+1}  \in R_{q+1} \Leftrightarrow  a^{\tau \mu}_k(x) f^k_{\mu +1_i}(x)+ { \partial}_i a^{\tau \mu}_k(x) f^k_{\mu}(x) =0 $ as identities on $X$ or at least over an open subset $U\subset X$. Differentiating the first relation with respect to $x^i$ and substracting the second, we finally obtain:  \\
\[  a^{\tau \mu}_k(x) ({\partial}_if^k_{\mu}(x) - f^k_{\mu +1_i}(x))=0 \Rightarrow df_{q+1}\in T^*\otimes R_q \]
and the Spencer operator restricts to $d:{\cal{R}}_{q+1} \rightarrow T^*\otimes R_q$. We set 
$ R^{(1)}_{q+r}={\pi}^{q+r+1}_{q+r}( R_{q+r+1})$.  \\

\noindent
{\bf DEFINITION 2.2}: The {\it symbol} of $R_q$ is the family $g_q=R_q \cap S_qT^*\otimes E$ of vector spaces over $X$. The symbol $g_{q+r}$ of ${\cal{R}}_{q+r}$ only depends on $g_q$ by a direct prolongation procedure. We may define the vector bundle $F_0$ over ${\cal{R}}_q$ by the short exact sequence $0 \rightarrow R_q \rightarrow J_q(E) \rightarrow F_0 \rightarrow 0$ and we have the exact induced sequence $0 \rightarrow g_q \rightarrow S_ qT^*\otimes E \rightarrow F_0$ .\\

When $\mid \mu \mid =q$, we obtain:  \\
 \[  g_q=\{ v^k_{\mu}\in S_qT^*\otimes E \mid a^{\tau \mu}_k(x)v^k_{\mu}=0\} , \mid \mu \mid = q   \]
  \[  \Rightarrow  g_{q+r}={\rho}_r(g_q)=\{ v^k_{\mu + \nu}\in S_{q+r} T^*\otimes E \mid a^{\tau \mu}_k(x)v^k_{\mu + \nu}=0\}, \mid \mu\mid=q, \mid \nu \mid =r  \]
In general, neither $g_q$ nor $g_{q+r}$ are vector bundles over $X$ as can be seen in the simple example $x y_x - y = 0 \Rightarrow x y_{xx}=0$.  \\

On ${\wedge}^sT^*$ we may introduce the usual bases $\{dx^I=dx^{i_1}\wedge ... \wedge dx^{i_s}\}$ where we have set 
$I=(i_1< ... <i_s)$. In a purely algebraic setting, one has:  \\

\noindent
{\bf PROPOSITION 2.3}: There exists a map $\delta:{\wedge}^sT^*\otimes S_{q+1}T^*\otimes E\rightarrow {\wedge}^{s+1}T^*\otimes S_qT^*\otimes E$ which restricts to $\delta:{\wedge}^sT^*\otimes g_{q+1}\rightarrow {\wedge}^{s+1}T^*\otimes g_q$ and ${\delta}^2=\delta\circ\delta=0$.\\

{\it Proof}: Let us introduce the family of s-forms $\omega=\{ {\omega}^k_{\mu}=v^k_{\mu,I}dx^I \}$ and set $(\delta\omega)^k_{\mu}=dx^i\wedge{\omega}^k_{\mu+1_i}$. We obtain at once $({\delta}^2\omega)^k_{\mu}=dx^i\wedge dx^j\wedge{\omega}^k_{\mu+1_i+1_j}=0$ and $ a^{\tau \mu}_k(\delta \omega)^k_{\mu}=dx^i \wedge(  a^{\tau \mu}_k{\omega}^k_{\mu +1_i})=0$.  \\
\hspace*{12cm} $\Box $  \\

The kernel of each $\delta$ in the first case is equal to the image of the preceding $\delta$ but this may no longer be true in the restricted case and we set:\\

\noindent
{\bf DEFINITION 2.4}: Let $B^s_{q+r}(g_q)\subseteq Z^s_{q+r}(g_q)$ and $H^s_{q+r}(g_q)=Z^s_{q+r}(g_q)/B^s_{q+r}(g_q)$ with $H^s(g_q)=H^s_q(g_q)$ be the coboundary space $im(\delta)$, cocycle space $ker(\delta)$ and cohomology space at ${\wedge}^sT^*\otimes g_{q+r}$ of the restricted $\delta$-sequence which only depend on $g_q$ and may not be vector bundles. The symbol $g_q$ is said to be s-{\it acyclic} if $H^1_{q+r}=...=H^s_{q+r}=0, \forall r\geq 0$, {\it involutive} if it is n-acyclic and {\it finite type} if $g_{q+r}=0$ becomes trivially involutive for r large enough. In particular, if $g_q$ is involutive {\it and} finite type, then $g_q=0$. Finally, $S_qT^*\otimes E$ is involutive for any $ q\geq 0$ if we set $S_0T^*\otimes E=E$. \\

Having in mind the example of $xy_x-y=0\Rightarrow xy_{xx}=0$ with rank changing at $x=0$, we have:  \\

\noindent
{\bf PROPOSITION 2.5}: If $g_q$ is $2$-acyclic and $g_{q+1}$ is a vector bundle, then $g_{q+r}$ is a vector bundle $ \forall r \geq1$. \\

\noindent
{\it Proof}: We may define the vector bundle $F_1$ by the following ker/coker exact sequence where we denote by $h_1\subseteq T^*\otimes F_0$ the image of the central map:  \\
\[   0 \rightarrow g_{q+1} \rightarrow S_{q+1}T^*\otimes E \rightarrow T^*\otimes F_0 \rightarrow F_1 \rightarrow 0  \]
and we obtain by induction on $r$ the following commutative and exact diagram of vector bundles:  \\
{\footnotesize \[  \begin{array}{lccccccc}
 & 0\hspace{3mm}& & 0 \hspace{3mm} && 0\hspace{3mm}& & 0 \hspace{3mm}  \\
 & \downarrow \hspace{3mm}& & \downarrow \hspace{3mm}& & \downarrow \hspace{3mm} & & \downarrow  \hspace{3mm}\\
0 \rightarrow & g_{q+r+1} & \rightarrow & S_{q+r+1}T^*\otimes E & \rightarrow  & S_{r+1}T^* \otimes F_0 & \rightarrow  & S_rT^*\otimes F_1  \\
  &  \downarrow \delta  & & \downarrow \delta & & \downarrow  \delta & & \downarrow \delta  \\
0 \rightarrow & T^*\otimes g_{q+r} & \rightarrow & T^* \otimes S_{q+r}T^* \otimes E & \rightarrow & T^*\otimes S_rT^*\otimes F_0 &\rightarrow & T^*\otimes S_{r-1}T^*\otimes F_1  \\
   &  \downarrow \delta  & & \downarrow \delta & & \downarrow  \delta & & \\
0 \rightarrow & {\wedge}^2T^*\otimes g_{q+r-1} & \rightarrow &{\wedge}^2T^*\otimes S_{q+r-1}T^* \otimes E & \rightarrow & {\wedge}^2T^*\otimes S_{r-1}T^*\otimes F_0 &  &   \\
  & \downarrow \delta &  & \downarrow \delta & & & &   \\
 & {\wedge}^3T^*\otimes S_{q+r-2} T^* \otimes E & = &   {\wedge}^3T^*\otimes S_{q+r-2} T^* \otimes E & & & &\\
  & & & & & & &
\end{array}   \] }  \\

A chase proves that the upper sequence is exact at $S_{r+1} T^* \otimes F_0$ whenever $g_q$ is $2$-acyclic by extending the diagram. The proposition finally follows by upper-semicontinuity from the relation:   \\
\[ dim(g_{q+r+1})+dim(h_{r+1})=m \, \, dim(S_{q+r+1}T^*)   \]
\hspace*{12cm}  $\Box$  \\   \\

\noindent
{\bf LEMMA 2.6}: If $g_q$ is involutive and $g_{q+1}$ is a vector bundle, then $ g_q$ is also a vector bundle. In this case, changing linearly the local coordinates if necessary, we may look at the maximum number $\beta$ of equations that can be solved with respect to $v^k_{n...n}$ and the intrinsic number $\alpha=m-\beta$ indicates the number of $y$ that can be given arbitrarily.  \\

Using the exactness of the preceding diagram and chasing in the following diagram:
{\footnotesize \[  \begin{array}{lcccccccc}
 && 0 & & 0  && 0 & & \\
 && \downarrow \ & & \downarrow & & \downarrow & &   \\
0 &\rightarrow & g_{q+r+1} & \rightarrow & S_{q+r+1}T^*\otimes E & \rightarrow  & S_{r+1}T^* \otimes F_0 & \rightarrow  & S_rT^*\otimes F_1  \\
 & &  \downarrow   & & \downarrow & & \downarrow   & &    \\
0 & \rightarrow & R_{q+r+1} & \rightarrow &  J_{q+r+1}(E ) & \rightarrow & J_{r+1}(F_0) &  &   \\
&  & \downarrow  &  & \downarrow  & & \downarrow & &  \\
  0 & \rightarrow & R_{q+r} & \rightarrow &  J_{q+r}(E )  &  \rightarrow & J_r(F_0) & &     \\
  & &  &  & \downarrow &  & \downarrow  & &    \\
 & & & &   0 & & 0  & &
\end{array}   \] }  \\
we have (See [24], p 95-98 for details): \\

\noindent
{\bf THEOREM 2.7}: If $ R_q\subset J_q(E)$ is a system of order $q$ on $ E $ such that $g_{q+1}$ is a vector bundle and $g_q$ is $2$-acyclic, then there is a commutative diagram: \\
\[  \begin{array}{rccccccl}
0 & \rightarrow & R^{(1)}_{q+r} & \rightarrow & R_{q+r} &  \stackrel{{\kappa}_r}{\longrightarrow} &  S_rT^*\otimes F_1  \\ 
   & & \downarrow & & \downarrow & & \downarrow &    \\
0 & \rightarrow &  J_r(R^{(1)}_q)  &  \rightarrow &  J_r(R_q)  &\stackrel{J_r(\kappa)}{ \longrightarrow}  & J_r(F_1) 
\end{array}  \]
where ${\kappa}_r$ is called the $r$-{\it curvature} and $\kappa={\kappa}_0$ is simply called the {\it curvature} of $ R_q$.  \\

We notice that $ R_{q+r+1}={\rho}_r (R_{q+1})$ and $ R_{q+r}=  {\rho}_r( R_q)$ in the following commutative diagram:  \\
\[  \begin{array}{ccc}
             R_{q+r+1}  &  \stackrel{{\pi}^{q+r+1}_{q+1}}{\longrightarrow }& R_{q+1}   \\
          \hspace{10mm}   \downarrow   {\pi}^{q+r+1}_{q+r}  &         &  \hspace{7mm}\downarrow {\pi}^{q+1}_q \\
             R^{(1)}_{q+r}    &  \stackrel{{\pi}^{q+r}_q }{\longrightarrow} & R^{(1)}_q    \\
               \cap   &   &   \cap   \\
               R_{q+r} &  \stackrel{{\pi}^{q+r}_q}{\longrightarrow}& R_q
               \end{array}  \]
We also have $ R^{(1)}_{q+r} \subseteq {\rho}_r( R^{(1)}_q)$ because we have successively:  \\
\[  \begin{array}{ccl}
R^{(1)}_{q+r}={\pi}^{q+r+1}_{q+r}(R_{q+r+1}) & = & {\pi}^{q+r+1}_{q+r}(J_r(R_{q+1}) \cap J_{q+r+1}(E)  \\
       & \subseteq  & J_r({\pi}^{q+1}_q) (J_r( R_{q+1})) \cap J_{q+r}(E)  \\
          & =  &  J_r(R^{(1)}_q) \cap J_{q+r}(E)  \\
              & = &  {\rho}_r(R^{(1)}_q)
              \end{array}   \]
while chasing in the following commutative $3$-dimensional diagram:  \\
\[  \begin{array}{rcccc}
  &  & J_r(R_{q+1}) & \longrightarrow & J_r(J_{q+1}(E))   \\
   & \nearrow  & \downarrow &  &  \nearrow \hspace{10mm} \\
   R_{q+r+1}  &   & \longrightarrow & J_{q+r+1}(E) &  \downarrow  \\
   \downarrow \hspace{5mm} &   &  J_r(R_q) & \longrightarrow  &  J_r(J_q(E))  \\
    & \nearrow &  & \downarrow & \nearrow  \hspace{10mm} \\
   R_{q+r} & & \longrightarrow & J_{q+r}(E)
   \end{array}  \]             
with a well defined map $J_r({\pi}^{q+1}_q):J_r(J_{q+1}(E)) \rightarrow J_r(J_q(E))$. \\

We finally obtain the following crucial Theorem (Compare to [21], p 72-74 or [27], p 340 to [14]):  \\
                             
\noindent
{\bf THEOREM 2.8}: Let $ R_q\subset J_q(E)$ be a system of order $q$ on $ E $ such that $ R_{q+1}$ is a vector sub-bundle of $J_{q+1}(E)$. If $g_q$ is $2$-acyclic and $g_{q+1}$ is a vector bundle, then $ R^{(1)}_{q+r}={\rho}_r (R^{(1)}_q), \forall r\geq 0$.  \\

\noindent
{\bf DEFINITION 2.9}: A system $ R_q\subset J_q(E)$ is said to be {\it formally integrable} if ${\pi}^{q+r+1}_{q+r}: R_{q+r+1}\rightarrow R_{q+r}$ is an epimorphism of vector bundles $\forall r\geq 1$ and {\it involutive} if it is formally integrable with an involutive symbol $g_q$. We have the following useful test [14, 21, 24, 27, 54]: \\ 

\noindent
{\bf COROLLARY 2.10}: Let $ R_q\subset J_q(E)$ be a system of order $q$ on $ E $ such that $ R_{q+1}$ is a vector sub-bundle of $J_{q+1}(E)$. If $g_q$ is $2$-acyclic (involutive) and if the map ${\pi}^{q+1}_q: R_{q+1} \rightarrow R_q$ is an epimorphism of vector bundles, then $ R_q$ is formally integrable (involutive). Such a result can be easily extended to nonlinear systems ([24]). \\

The next procedure providing a {\it Pommaret basis} and where  {\it one may have to change linearly the independent variables if necessary}, is intrinsic even though it must be checked in a particular coordinate system called $\delta$-{\it regular} ([21]).  \\

$\bullet$ {\it Equations of class} $n$: Solve the maximum number ${\beta}^n_q$ of equations with respect to the jets of order $q$ and class $n$. Then call $(x^1,...,x^n)$ {\it multiplicative variables}.\\

$\bullet$ {\it Equations of class} $i\geq 1$: Solve the maximum number ${\beta}^i_q$ of {\it remaining} equations with respect to the jets of order $q$ and class $i$. Then call $(x^1,...,x^i)$ {\it multiplicative variables} and $(x^{i+1},...,x^n)$ {\it non-multiplicative variables}.\\

$\bullet$ {\it Remaining equations equations of order} $\leq q-1$: Call $(x^1,...,x^n)$ {\it non-multiplicative variables}.\\

\noindent
In actual practice, we shall use a {\it Janet tabular} where the multiplicative "variables" are in upper left position while the non-multiplicative variables are represented by dots in lower right position. ccording to the previous results,  a system of PD equations is {\it involutive} if its first prolongation can be obtained by prolonging its equations only with respect to the corresponding multiplicative variables. In that case, we may introduce the {\it characters} ${\alpha}^i_q=m\frac{(q+n-i-1)!}{(q-1)!((n-i)!}-{\beta}^i_q$ for $i=1, ..., n$ with ${\alpha}^1_q\geq ... \geq {\alpha}^n_q\geq 0$ and we have $dim(g_q)={\alpha}^1_q+...+ {\alpha}^n_q$ while $dim(g_{q+1})={\alpha}^1_q+...+ n {\alpha}^n_q$. \\
 
\noindent 
We now recall the main results and definitions that are absolutely needed for the applications.\\
With canonical epimorphism ${\Phi}_0=\Phi:J_q(E) \Rightarrow J_q(E)/R_q=F=F_0$, the various prolongations are described by the following commutative and 
exact "{\it introductory diagram} " in which we set $R'_r =im({\rho}_r(\Phi))  \subset J_r(F_0)$ with $R'_0 = F_0$, $Q_r = coker({\rho}_r(\Phi))$ and 
$h_{r+1} = coker ({\sigma}_{r+1}(\Phi))$: \\
\[ \fbox{ $   \begin{array}{rcccccccccccl}
   & &  0  &  &  0  &  &  0  &  &  0  & & & \\
   &  &  \downarrow &  & \downarrow &  &  \downarrow  &  &  \downarrow    &  &  & \\ 
0  & \rightarrow  &  g_{q+r+1}  & \rightarrow & S_{q+r+1}T^*\otimes E & \stackrel{{\sigma}_{r+1}(\Phi)}{\rightarrow} & g'_{r+1}  & \subset &   S_{r+1}T^*\otimes F_0 & \rightarrow  & h_{r+1} & \rightarrow &  0  \\
&  & \downarrow &  & \downarrow &  &  \downarrow  &  &  \downarrow   & & \downarrow &     \\
0 & \rightarrow &  R_{q+r+1} & \rightarrow  & J_{q+r+1}(E) & \stackrel{{\rho}_{r+1}(\Phi)}{\rightarrow}  & R'_{r+1}&  \subset  &   J_{r+1}(F_0)  & \rightarrow & Q_{r+1} & \rightarrow & 0 \\
 &   &  \downarrow &  & \downarrow &  &  \downarrow  &  & \downarrow    & & \downarrow &   \\
0 &  \rightarrow &  R_{q+r} & \rightarrow  & J_{q+r}(E) & \stackrel{{\rho}_r(\Phi)}{\rightarrow}  & R'_r  & \subset &   J_r(F_0)  & \rightarrow & Q_r & \rightarrow &  0  \\ 
     &  & &  &  \downarrow  &  & \downarrow  &  &  \downarrow   &  &  \downarrow &  \\
     &  &  & & 0  & & 0  &  & 0 & & 0 &
\end{array}  $  } \]
Chasing along the diagonal of this diagram while applying the standard "{\it snake}" lemma, we notice that $im({\sigma}_r(\Phi)) \subseteq g'_r$ and obtain the useful "{\it long exact connecting sequence} ":  \\
 \[  0  \rightarrow  g_{q+r+1}  \rightarrow R_{q+r+1}  \rightarrow R_{q+r}  \rightarrow h_{r+1} \rightarrow Q_{r+1}  \rightarrow Q_r \rightarrow 0  \]
which is thus connecting in a tricky way FI ({\it lower left}) with CC ({\it upper right}).  \\
\noindent
A key step in the procedure for constructing differential sequences will be to use the following theorems and corollary (See [21, 27] for acyclicity or involutivity):  \\

\noindent
{\bf THEOREM 2.11}: There is a finite {\it Prolongation/Projection} (PP) algorithm providing two integers $r,s\geq 0$ by successive increase of each of them such that the new system $R^{(s)}_{q+r}= {\pi}^{q+r+s}_{q+r}(R_{q+r+s})$ has the same solutions as $R_q$ but is FI with a $2$-acyclic or involutive symbol and first order CC. The maximum order of ${\cal{D}}_1$ is thus equal to $ r+s+1$ as we used $r+s$ prolongations but it may be lower because certain CC may generate the higher order ones as will be seen in the motivating examples. As long as this procedure has not been achieved, {\it nothing} can be said about the CC (Fine examples can be found in [37] and the recent [45]). \\

\noindent 
{\bf DEFINITION 2.12}: A differential sequence is said to be {\it formally exact} if it is exact on the jet level composition of the prolongations involved. A formally exact sequence is said to be {\it strictly exact} if all the operators/systems involved are FI (See [34] for more details). A strictly exact sequence is called {\it canonical} if all the operators/systems are involutive. \\

When $d: J_{q+1}(E) \rightarrow T^* \otimes J_q(E): f_{q+1} \rightarrow j_1(f_q)  - f_{q+1}$ is the Spencer operator, we have:  \\

\noindent
{\bf PROPOSITION 2.13}: If $R_q\subset J_q(E)$ and $R_{q+1} \subset J_{q+1}(E)$ are two systems of respective orders $q$ and $q+1$, then $R_{q+1} \subset {\rho}_1(R_q)$ if and onlty if ${\pi}^{q+1}_q (R_{q+1}) \subset R_q$ {\it and} $dR_{q+1} \subset T^* \otimes R_q$.  \\

\noindent
{\bf DEFINITION 2.14}: Let us "cut" the preceding introductory diagram  by means of a central vertical line and define $R'_r = im({\rho}_r(\Phi)) \subseteq J_r(F_0)$ with $R'_0=F_0$. Chasing in this diagram, we notice that ${\pi}^{r+1}_r:J_{r+1}(F_0) \rightarrow J_r(F_0)$ induces an epimorphism ${\pi}^{r+1}_r: R'_{r+1} \rightarrow R'_r, \forall r\geq 0$. However, a chase in this diagram proves that {\it the kernel of this epimorphism is not} $im({\sigma}_{r+1}(\Phi))$ unless $R_q$ is FI ({\it care}). For this reason, we shall define it to be {\it exactly} 
$g'_{r+1}$.  \\

\noindent
{\bf THEOREM 2.15}:  $R'_{r+1} \subseteq {\rho}_1(R'_r)$ and $ dim ({\rho}_1(R'_r)) - dim (R'_{r+1})$ is the number of new generating CC of order $r+1$ .  \\

\noindent
{\bf COROLLARY; 2.16}: The system $R'_r \subset J_r(F_0)$ becomes FI with a $2$-acyclic or involutive symbol and $R'_{r+1}={\rho}_1(R'_r)  \subset J_{r+1}(F_0)$ when $r$ is large enough.  \\ 

Using the preceding intrinsic results, one may always suppose that we may start with an involutive system $R_q \subset  J_q(E)$, that is formaly integrable with an involutive symbol $g_q \subset S_q T^* \otimes E$. Then, using the Spencer operator, one can construct {\it STEP BY STEP} or {\it AS A WHOLE}, 
the three differential sequences related by the following {\it fundamental diagram I} below, namely the Specer sequenc, the Janet sequence and the central {\it hybrid sequence} which is at the same time the Janet sequence for the trivially involutive operator $j_q$ and the Spencer sequence for the first order system $J_{q+1}(E) \subset J_1 ( J_q(E))$. With more details, when $R_q$ is involutive, the operator ${\cal{D}}:E\stackrel{j_q}{\rightarrow} J_q(E)\stackrel{\Phi}{\rightarrow} J_q(E)/R_q=F_0$ of order $q$ is said to be {\it involutive}. Introducing the {\it Janet bundles} $F_r= {\wedge}^rT^*\otimes J_q(E)/({\wedge}^rT^*\otimes R_q + \delta (S_{q+1}T^*\otimes E))$, we obtain the {\it linear Janet sequence} induced by the Spencer operator that has been introduced in ([21, 24]):\\
\[  0 \longrightarrow  \Theta \longrightarrow E \stackrel{\cal{D}}{\longrightarrow} F_0 \stackrel{{\cal{D}}_1}{\longrightarrow}F_1 \stackrel{{\cal{D}}_2}{\longrightarrow} ... \stackrel{{\cal{D}}_n}{\longrightarrow} F_n \longrightarrow 0   \]
where each other operator is first order involutive and generates the CC of the preceding one. \\
Similarly, introducing the {\it Spencer bundles} $C_r= {\wedge}^rT^* \otimes R_q / \delta ({\wedge}^{r-1}T^* \otimes g_{q+1})$ we obtain the {\it linear Spencer sequence} induced by the Spencer operator [21, 24, 54]:  \\
\[   0\longrightarrow \Theta \stackrel{j_q}{\longrightarrow} C_0 \stackrel{D_1}{\longrightarrow} C_1 \stackrel{D_2}{\longrightarrow}... \stackrel{D_n}{\longrightarrow} C_n 
\longrightarrow  0  \]   
The Janet and Spencer bundles are respectively defined by the formulas:  
\[ \fbox{  $  0 \rightarrow {\wedge}^rT^*\otimes R_q + \delta({\wedge}^{r-1}T^*\otimes S_{q+1}T^*\otimes E ) \rightarrow  {\wedge}^rT^*\otimes J_q(E) \rightarrow  F_r \rightarrow 0 $ } \]
\[  \fbox{  $   0 \rightarrow \delta({\wedge}^{r-1}T^* \otimes g_{q+1}) \rightarrow {\wedge} ^rT^*\otimes R_q \rightarrow C_r  \rightarrow 0  $  } \]
The two sequences are related by the following commutative diagram with short exact vertical sequences:  \\

  \[  \begin{array}{rccccccccccccccl}
 &&&&& 0 &&0&&0&  & 0 & & \\
 &&&&& \downarrow && \downarrow && \downarrow  & & \downarrow  &  &\\
  & 0& \rightarrow& \Theta &\stackrel{j_q}{\rightarrow}& C_0  &\stackrel{D_1}{\rightarrow}&  C_1 &\stackrel{D_2}{\rightarrow} ... \stackrel{D_{n-1}}{\rightarrow}&  C_{n-1} &\stackrel{D_n}{\rightarrow}  & C_n & \rightarrow  & 0  \\
  &&&&& \downarrow & & \downarrow & & \downarrow & & \downarrow &   & \\
   & 0 & \rightarrow &  E & \stackrel{j_q}{\rightarrow} & C_0(E)  & \stackrel{D_1}{\rightarrow} &  C_1(E)  &\stackrel{D_2}{\rightarrow}... \stackrel{D_{n-1}}{\rightarrow} & C_{n-1}(E) &  \stackrel{D_n}{\rightarrow} & C_n(E) & \rightarrow  & 0 \\
   & & & \parallel && \hspace{5mm}\downarrow {\Phi}_0 & &\hspace{5mm} \downarrow {\Phi}_1 & & \hspace{5mm}\downarrow {\Phi}_{n-1} &  &\hspace{5mm}\downarrow {\Phi}_n & \\
   0 \rightarrow & \Theta &\rightarrow &  E & \stackrel{\cal{D}}{\rightarrow} &  F_0  & \stackrel{{\cal{D}}_1}{\rightarrow} &  F_1  & \stackrel{{\cal{D}}_2}{\rightarrow} ... \stackrel{{\cal{D}}_{n-1}}{\rightarrow}&  F_{n-1} & \stackrel{{\cal{D}}_n}{\rightarrow}&  F_n  & \rightarrow  &  0      \\
   &&&&& \downarrow & & \downarrow & & \downarrow &  & \downarrow &    &   \\
   &&&&& 0 && 0 && 0  &  &  0  &  &
   \end{array}     \]
We finally recall that, when $R_q \subset J_q(E)$ is an involutive system, then the Spencer sequence is nothing else than the Janet sequence for the first order system $R_{q+1} \subset J_1(R_q)$. Such a (difficult) result can be obtained by using inductively a snake chase in the following commutative and exact diagram, starting with ${\Phi}_0 = \Phi$ with $R_q \subset J_q(E)$ for $r=0$ and we have for $r=1$ the commutative and exact diagram allowing to define $C_1$:
\noindent
\[   \begin{array}{rcccccccl}
  && 0 && 0 &&  0 & &  \\
&  &\downarrow  &  &\downarrow  &  & \downarrow  &  &  \\
0 & \rightarrow & g_{q+1} & \rightarrow & T^* \otimes R_q & \rightarrow & C_1 & \rightarrow & 0  \\
    &  & \downarrow  &  & \downarrow  &  &  \parallel &  &  \\
0 & \rightarrow &   R_{q+1}  &  \rightarrow & J_1(R_q) &  \rightarrow & C_1 & \rightarrow & 0 \\
&  & \downarrow  &  &\downarrow   &  & \downarrow   &  &\\
0 & \rightarrow &  R_q & = & R_q & \rightarrow & 0 & & \\
   &  & \downarrow  &  & \downarrow  &  &   &  & \\
   &  & 0 &  & 0 &  &  &  &  
\end{array}  \]
More generally, one has $R_{q+r} \subset J_r(R_q)$ for $r \geq 1$ and the procedure is ending when $r=n$ because all the $\delta$-sequences are exact at ${\wedge}^s T¨^* \otimes g_{q+r}$ for any $ 0 \leq s \leq n, r\geq 0$:
\[ \scriptsize { \begin{array}{rccccccccccccl}
   & &  &  & &  0  &  & 0  &   &  &  & 0  &  &   \\
   & &  & & & \downarrow & & \downarrow & &  &  & \downarrow  &  &   \\
&  0 & \rightarrow &  R_{q+r} & \rightarrow & J_r(R_q) &\rightarrow & J_{r-1}(C_1) &\rightarrow & ...  &\rightarrow   & C_r & \rightarrow  0 \\
  & &  &   & & \downarrow & & \downarrow &  &  &  & \downarrow  &  \\
   & 0 & \rightarrow &  J_{q+r}(E) & \rightarrow & J_r(J_q(E)) & \rightarrow  &  J_{r-1}(C_1(E))  & \rightarrow  &... &\rightarrow   & C_r(E) & \rightarrow  0  \\
  & &  &  \parallel & & \hspace{10mm} \downarrow J_r({\Phi}_0) & & \hspace{13mm} \downarrow J_{r-1}({\Phi}_1)& & &  & \hspace{5mm} \downarrow  {\Phi}_r     \\
0  \rightarrow & R_{q+r} & \rightarrow &J_{q+r}(E) & \rightarrow & J_r(F_0) &\rightarrow & J_{r-1}(F_1) & \rightarrow &  ... & \rightarrow & F_r & \rightarrow  0  \\
 & &  &  \downarrow & & \downarrow  & & \downarrow  & & &  & \downarrow  &   \\
   & &  & 0  & &  0  &  & 0  &   &   & &  0 &   
\end{array} } \]   \\

All the following motivating example are taken from standard papers comparing Janet bases and Pommaret  bases but their applications need no comment. Each one is illustrating a specific singular result that cannot be obtained by classical techniques  ([2, 3, 5, 8-13, 51, 56]).  \\

\noindent
{\bf EXAMPLE  2.17 }: (Proceedings GIFT 2006, p 187) With polynomial variables $(x,y,z)$ and ground field $k = \mathbb{Q}$, let us consider the ideal $ \mathfrak{a} = (z^2  - y^2 - 2 x^2, x z + x y, y z + y^2 + x^2) \subset k[x, y, z] = A$ with ordering $x \prec y \prec z$, which is not prime because $x(y + z) \in \mathfrak{a}$. Also, adding twice the third polynomial to the first, we obtain $(y + z)^2 \in \mathfrak{a}$ and the prime ideal $rad(\mathfrak{a}) = \sqrt{\mathfrak{a}} = (y + z, x) = \mathfrak{p}$. Now, passing to jet notations in the linear PD framework, we get the following linear homogeneous system of order two $R_2 \subset J_2 (E)$ with $dim(R_2) = dim(J_2(E)) - 3 = 7$: \\
\[  \left\{\begin{array}{lcl}
y_{33} - y_{22}  - 2 y_{11} & =   &  0  \\
y_{23} + y_{22} + y_{11} & =  &  0 \\
y_{13} + y_{12} & = & 0
\end{array} 
\fbox{ $  \begin{array}{ccc}
1 & 2 & 3 \\
1 & 2 & \bullet  \\
1 & \bullet & \bullet
\end{array} $ }   \right. \]
Using the last bottom two non-multiplicative variables, we let the reader check that {\it the system is not involutive in this system of coordinates }, which does not mean it is not involutive by choosing convenient $\delta$-regular coordinates, even if we already know that $R_2$ is formally integrable (FI) because it is homogeneous of order two. One cannot proceed ahead for computing the characters, even if we know that ${\alpha}^3_2=0$ because the differential module admits the torsion element $z=y_3 + y_2$ with $d_1 z=0$ and the first equation is solved with respect to $y_{33}$. Looking at the previous prime ideal, we may add twice the second equation to the first and choose the new variables:
\[  {\bar{x}}^3 = x^3 + x^2, \,\,\, {\bar{x}}^2 = x^1, \,\,\, {\bar{x}}^1= x^2  \]
in order to obtain now the new equivalent following system $R_2 \subset J_2(E)$ after taking out the bar for simplicity (compare to GIFT):
\[  \left\{\begin{array}{lcl}
y_{33}  & =   &  0  \\
y_{23} & =  &  0 \\
y_{22} + y_{13} & = & 0
\end{array} 
\fbox{ $  \begin{array}{ccc}
1 & 2 & 3 \\
1 & 2 & \bullet  \\
1 & 2& \bullet
\end{array} $ }   \right. \]
As $dim(E) = 1$, its symbol $g_2 $ is involutive with $dim(g_2)= 3$ with parametric jets $(y_{11}, y_{12}, y_{13})$. The three characters are thus ${\alpha}^3_2=0, {\alpha}^2_2=0, {\alpha}^1_2= 3$ and we check that $dim(g_2) = {\alpha}^1_2 + {\alpha}^2_2 + {\alpha}^3_2 =3$. We have also $dim(g_3)= {\alpha}^1_2 + 2 {\alpha}^2_2 + 3 {\alpha}^3_2= {\alpha}^1_2 =3$ and so on with $dim(g_{2+r})=3, \forall r\geq 0$. \\
As $R_2$ is an involutive system, the corresponding formally exact Janet sequence can be written as follows with only two differentially independent CC for $\cal{D}$:
\[     0 \longrightarrow  \Theta \longrightarrow E \underset 2{\stackrel{{\cal{D}}}{\longrightarrow} }F_0 \underset 1{\stackrel{{\cal{D}}_1}{\longrightarrow}} F_1 \longrightarrow 0  \]
and thus $dim(E)= 1, dim(F_0) = 3, dim(F_1) =2$.  \\
For a later use, we write down the following useful Spencer $\delta$-sequences  which are exact because $g_2$ is involutive (See [21] , Chapter 3, Sections 1 and 2 for details) while providing the dimensions:
\[     0 \longrightarrow  g_3 \stackrel{\delta}{\longrightarrow} T^* \otimes g_2 \stackrel{\delta}{\longrightarrow} {\wedge}^2 T^* \otimes T^* \otimes E \Rightarrow
0 \longrightarrow  3 \stackrel{\delta}{\longrightarrow} 9 \stackrel{\delta}{\longrightarrow} 9  \]
\[       0 \longrightarrow  g_4 \stackrel{\delta}{\longrightarrow} T^* \otimes g_3 \stackrel{\delta}{\longrightarrow} {\wedge}^2 T^* \otimes g_2 \stackrel{\delta}{\longrightarrow} {\wedge}^3 T^* \otimes T^* \otimes E  \Rightarrow    0 \longrightarrow  3 \stackrel{\delta}{\longrightarrow} 9 \stackrel{\delta}{\longrightarrow} 9 \stackrel{\delta}{\longrightarrow}  3 \longrightarrow 0 \]
 \[ 0 \longrightarrow  g_5 \stackrel{\delta}{\longrightarrow} T^* \otimes g_4 \stackrel{\delta}{\longrightarrow}  {\wedge}^2 T^* \otimes g_3 \stackrel{\delta}{\longrightarrow} {\wedge}^3 T^* \otimes g_2 \longrightarrow 0 \Rightarrow 0 \longrightarrow 3 \stackrel{\delta}{\longrightarrow} 9 \stackrel{\delta}{\longrightarrow}  9 \stackrel{\delta}{\longrightarrow} 3 \longrightarrow 0  \]
\[       0 \longrightarrow  g_6 \stackrel{\delta}{\longrightarrow} T^* \otimes g_5 \stackrel{\delta}{\longrightarrow} {\wedge}^2 T^* \otimes g_4 \stackrel{\delta}{\longrightarrow} {\wedge}^3 T^* \otimes g_3  \longrightarrow 0  \Rightarrow 
  0 \longrightarrow  3 \stackrel{\delta}{\longrightarrow} 9 \stackrel{\delta}{\longrightarrow}9 \stackrel{\delta}{\longrightarrow} 3  \longrightarrow 0 \]
We obtain therefore easily the Spencer bundles $C_r = {\wedge}^r T^* \otimes R_2 / \delta ({\wedge}^{r-1} T^* \otimes g_3)$. Hence, taking into account the exactness of the previous $\delta$- sequences, we get successively: 
\[  C_0 = R_2, C_1 =T^* \otimes R_2/\delta(g_3), C_2 = {\wedge}^2T^* \otimes R_2 / \delta ( T^* \otimes g_3), C_3 = {\wedge}^3 T^* \otimes R_2 / \delta ({\wedge}^2 T^* \otimes g_3)  \]
and thus $dim(C_0)= 7, dim(C_1)= (3 \times 7) - 3 = 18, dim(C_2) = (3 \times 7) - 6 = 15, dim(C_3) = 7 - 3 = 4 $.  \\

 \[  \begin{array}{rccccccccccccccl}
 &&&&& 0 &&0&&0&  & 0 & & \\
 &&&&& \downarrow && \downarrow && \downarrow  & & \downarrow  &  &\\
  & 0& \longrightarrow& \Theta &\stackrel{j_2}{\longrightarrow}& C_0  &\stackrel{D_1}{\longrightarrow}&  C_1 &\stackrel{D_2}{\longrightarrow} &  C_2 &\stackrel{D_3}{\longrightarrow}  & C_3 & \longrightarrow  & 0  \\
  &&&&& \downarrow & & \downarrow & & \downarrow & & \downarrow &   & \\
   & 0 & \longrightarrow &  E & \stackrel{j_2}{\longrightarrow} & C_0(E)  & \stackrel{D_1}{\longrightarrow} &  C_1(E)  &\stackrel{D_2}{\longrightarrow} & C_2(E) &  \stackrel{D_3}{\longrightarrow} & C_3(E) & \longrightarrow  & 0 \\
   & & & \parallel && \hspace{5mm}\downarrow {\Phi}_0 & &\hspace{5mm} \downarrow {\Phi}_1 & & \hspace{5mm}\downarrow {\Phi}_2 &  &\hspace{5mm}\downarrow {\Phi}_3 & \\
   0 \longrightarrow & \Theta &\longrightarrow &  E & \stackrel{\cal{D}}{\longrightarrow} &  F_0  & \stackrel{{\cal{D}}_1}{\longrightarrow} &  F_1  & \stackrel{{\cal{D}}_2}{\longrightarrow} &  F_2 & \stackrel{{\cal{D}}_3}{\longrightarrow}&  F_3  & \longrightarrow  &  0      \\
   &&&&& \downarrow & & \downarrow & & \downarrow &  & \downarrow &    &   \\
   &&&&& 0 && 0 && 0  &  &  0  &  &
   \end{array}     \]
  \[  \begin{array}{rccccccccccccccl}
 &&&&& 0 &&0&&0&  & 0 & & \\
 &&&&& \downarrow && \downarrow && \downarrow  & & \downarrow  &  &\\
  & 0& \longrightarrow& \Theta &\stackrel{j_2}{\longrightarrow}& 7 &\stackrel{D_1}{\longrightarrow}& 18 &\stackrel{D_2}{\longrightarrow} &  15 &\stackrel{D_3}{\longrightarrow}  & 4 & \longrightarrow  & 0  \\
  &&&&& \downarrow & & \downarrow & & \downarrow & & \downarrow &   & \\
   & 0 & \longrightarrow &  1 & \stackrel{j_2}{\longrightarrow} & 10  & \stackrel{D_1}{\longrightarrow} & 
   20  &\stackrel{D_2}{\longrightarrow} & 15 &  \stackrel{D_3}{\longrightarrow} & 4 & \longrightarrow  & 0 \\
   & & & \parallel && \hspace{5mm}\downarrow {\Phi}_0 & &\hspace{5mm} \downarrow {\Phi}_1 & &  \downarrow  &  &  \downarrow  & \\
   0 \longrightarrow & \Theta &\longrightarrow &  1 & \stackrel{\cal{D}}{\longrightarrow} &  3  & \stackrel{{\cal{D}}_1}{\longrightarrow} &  2  & \longrightarrow &  0 
   & &  0  & &       \\
   &&&&& \downarrow & & \downarrow & & &  & &    &   \\
   &&&&& 0 && 0 &&   &  &    &  &
   \end{array}     \] 
The morphisms ${\Phi}_1, {\Phi}_2, {\Phi}_3$ in the vertical short exact sequences are inductively induced from the morphism ${\Phi}_0= \Phi$ in the first short exact vertical sequence on the left. The central horizontal sequence can be called " {\it hybrid sequence} " because it is at the same time a Spencer sequence for the first order system $J_3(E) \subset J_1(J_2(E))$ over $J_2(E)$ and a formally exact Janet sequence for the involutive injective operator $j_2:E \rightarrow J_2(E)$. It can be constructed step by step, starting with the short exact sequence:   $0 \longrightarrow J_3(E) \longrightarrow J_1(J_2(T)) \longrightarrow C_1(E) \longrightarrow 0$ or, equivalently, the short exact symbol sequence:  $0 \longrightarrow 
S_3 T^* \otimes E  \longrightarrow T ^* \otimes J_2 (E) \longrightarrow C_1(E) \longrightarrow 0$, kernel of the projection onto $J_2(T)$ of these affine vector bundles. We also invite the reader, as an exercise, to construct it as a whole by introducing the Spencer bundles $ C_r(E) = {\wedge}^r T^* \otimes J_2(E) / \delta ({\wedge}^{r-1} T^* \otimes S_3 T^* \otimes E)$. The comparison with ([GIFT]) needs no comment on the usefulness and effectiveness of constructing both differential sequences. {\it In this particular case}, the Janet sequence may be quite simpler than the Spencer sequence, contrary to what could happen in other examples. We recall that they are {\it absolutely} needed for studying the group of conformal transformations as in ([43, 44]).  \\

\noindent
{\bf EXAMPLE 2.18}: ([13]) With polynomial variables $(x,y,z)$ and ground field $k = \mathbb{Q}$, let us consider the ideal $ \mathfrak{a} = (xy, xz, yz) \Rightarrow 
{\mathfrak{a}} = (x,y) \cap (x, z) \cap (y, z) = {\mathfrak{p}}_1 \cap {\mathfrak{p}}_2 \cap {\mathfrak{p}}_3$. Passing to jet notations in the linear PD framework, we get the following linear homogeneous system of order two $R_2 \subset J_2 (E)$ with $dim(R_2) = dim(J_2(E)) - 3 = 10 - 3 = 7$ defined by $y_{12}=0, y_{13}=0, y_{23}=0$ which is surely formally integrable as it is homogeneous. The symbol may not be involutive because there is no equation of class $3$, one equation of class $2$ and two equations of class $1$. Like in ([21]), let us make the change of coordinates $({\bar{x}}^1 = x^1, {\bar{x}}^2 = x^1 + x^2, {\bar{x}}^3 = x^2 + x^3) \Rightarrow (d_1 \rightarrow d_1 + d_2, d_2 \rightarrow d_2 + d_3, d_3 \rightarrow d_3) $. {\it Surprisingly}, dropping the bar, we get the following involutive system:  \\
 \[  \left\{  \begin{array}{lcl}
y_{33}  -  y_{13}  & =   &  0  \\
y_{23}  + y_{13}  & =  &  0 \\
y_{22}  + y_{12}  & = & 0  
\end{array} 
\fbox{ $  \begin{array}{ccc}
1 & 2 & 3 \\
1 & 2 & \bullet  \\
1 & 2 & \bullet  
\end{array} $ }   \right. \]
a result leading to the characters ${\alpha}^3_2 = 0, {\alpha}^2_2 = 0, {\alpha}^1_2 = 3 - 0 = 3$ and the Janet sequence:
\[   0 \rightarrow \Theta \rightarrow E \underset 2{\stackrel{{\cal{D}}}{\longrightarrow}}  F_1  \underset 1{\stackrel{{\cal{D}}}{\longrightarrow}} F_2 \rightarrow 0 \]
with Euler-Poincar\'{e} characteristic $rk_D(M) = 1 - 3 + 2 = 0$ in a coherent way. The interest of this example is that we obtain {\it exactly} the same formal results as in the preceding example though we have now a perfect ideal which is equal to its radical. We let the reader prove directly that $R_3 \subset J_1(R_2)$ is a first order involutive system. In actual practice, once we know the system is involutive, we need not change the coordinates backwards. However, we just need to make the change of variables 
$(z^1=y, z^2 = y_1, z^3 = y_2, z^4 = y_3) $ and obtain the new involutive system:
 \[  \left\{  \begin{array}{l}
z^1_3 - z^4 =0, z^2_3 - z^4_1 = 0, z^3_3 + z^4_1=0, z^4_3 - z^4_1 = 0  \\
z^1_2 - z^3 = 0, z^2_2 - z^3_1 = 0, z^3_2 + z^3_1 =0, z^4_2 + z^4_1 = 0 \\
z^1_1 - z^2 = 0
\end{array} 
\fbox{ $  \begin{array}{ccc}
1 & 2 & 3 \\
1 & 2 & \bullet  \\
1 & \bullet & \bullet  
\end{array} $ }   \right. \] 
We deduce that the module $ A / \mathfrak{a} $ is $2$-pure but this is out of the scope of this paper (See [30] for details). Contrary to what one could believe, this later system cannot be used in order to construct the Spencer sequence  (See [34] , Remark 2.3). This example is one of the best we know showing that, when one knows that a system is involutive {\it in one coordinate system}, one needs not go back to {\it this} coordinate system in order to construct sequences by means of diagram chasing. The case of Maxwell parametrization in elasticity will be similar later on. \\  

\noindent
{\bf EXAMPLE  2.19}: ([2]) With polynomial variables $(x,y,z)$ and ground field $k = \mathbb{Q}$, let us consider the ideal $ \mathfrak{a} = (z^2, y^2, z+x) \Rightarrow \sqrt{\mathfrak{a}} = (x, y, z) = \mathfrak{m}$ a maximal ideal. Dropping the present notation while passing to jet notations in the linear PDE framework, we get the following linear homogeneous system of order two $R_2 \subset J_2 (E)$ with $dim(R_2) = dim(J_2(E)) - 3 = 10 - 3 = 7$ defined by $y_{33}=0, y_{22}=0, y_3 + y_1=0$ which is not involutive because it is neither formally integrable and the symbol $g_2$ defined by $y_{33}=0, y_{22}=0$ with $y=0, y_1=0, y_2= 0$ is surely not involutive. We may thus use at least one prolongation, with the hope that the symbol $g_3$ of $R_3 = {\rho}_1(R_2)$ becomes involutive. We have the Janet tabular for $g_3$:
 \[  \left\{  \begin{array}{lcl}
y_{333}  & =   &  0  \\
y_{233}  & =  &  0 \\
y_{223} & = & 0  \\
y_{222} & = & 0  \\
y_{133}  &  = & 0 \\
y_{122}  & = & 0  
\end{array} 
\fbox{ $  \begin{array}{ccc}
1 & 2 & 3 \\
1 & 2 & \bullet  \\
1 & 2 & \bullet  \\
1 & 2 & \bullet \\
1 & \bullet & \bullet  \\
1 & \bullet & \bullet
\end{array} $ }   \right. \]
Using the non-multiplicative variables, we let the reader check that $g_3$ {\it  is indeed involutive in this system of coordinates } which is thus $\delta$-regular. We already know that $R_2$ and thus $R_3$ are surely {\it not} formally integrable (FI) because of the only PD equation $y_3 + y_1 = 0$. Using the important Theorems 2.8, 2.9, we are sure that $ {\rho}_r(R^{(1)}_3)= R^{(1)}_{r+3}, \forall r \geq 0$. However, one cannot proceed ahead for computing the characters, even if we know that ${\alpha}^3_2=0$ because the first equation is solved with respect to $y_{33}$. We have ${\alpha}^3_3=0, {\alpha}^2_3= 0, {\alpha}^1_3= 6 - 2=4$ because we have $6$ jet coordinates of class $1$, namely $(y_{111}, y_{112}, y_{113}, y_{122}, y_{123}, y_{133})$. We may thus construct $R^{(1)}_3$ with $dim(R^{(1)}_3)= 20 - 16 = 4  $ with parametric jets $(y, y_1, y_2, y_{12})$:  \\
 \[  \left\{  \begin{array}{lcl}
y_{333}  & =   &  0  \\
y_{233}  & =  &  0 \\
y_{223} & = & 0  \\
y_{222} & = & 0  \\
y_{133}  &  = & 0 \\
y_{123}  & = & 0  \\
y_{122} & = &  0 \\
y_{113} & = & 0 \\
y_{112} & = & 0 \\
y_{111} & =  & 0  \\
y_{33}  & = & 0  \\
y_{23} + y_{12} & = & 0 \\
y_{22} & = & 0 \\
y_{13} & = & 0  \\
y_{11} & = & 0 \\
y_3  + y_1 & = & 0  
\end{array} 
\fbox{ $  \begin{array}{ccc}
1 & 2 & 3 \\
1 & 2 & \bullet  \\
1 & 2 & \bullet  \\
1 & 2 & \bullet \\
1 & \bullet & \bullet  \\
1 & \bullet & \bullet  \\
1 & \bullet & \bullet  \\
1 & \bullet & \bullet \\
1 & \bullet & \bullet \\
1 & \bullet & \bullet  \\
\bullet & \bullet & \bullet \\
\bullet & \bullet & \bullet \\
\bullet & \bullet & \bullet \\
\bullet & \bullet & \bullet \\
\bullet & \bullet & \bullet \\
\bullet & \bullet & \bullet 
\end{array} $ }   \right. \]
{\it Strikingly}, we discover that the symbol of $R^{(1)}_3 $ is $g^{(1)}_3=0$, providing a finite type involutive system. Such a result could have been found directly through an explicit integration of the initial system which is providing the general solution in the form $y= a (x^2 x^3 - x^1 x^2) + bx^2 + c (x^3 - x^1)  + d$ with $4$ arbitrary constants $(a, b, c, d)$. 
We are thus in position to exhibit the two corresponding Janet and Spencer differential sequences, {\it without any reference to other technical tools} like in ([5]). The comparison needs no comment because, {\it in this particular case}, the Spencer sequence is isomorphic to the tensor product of the Poincar\'{e} differential sequence $(grad, curl, div)$ by 
${\mathbb{R}}^4$, a result {\it not evident at first sight}, while the Janet sequence is quite more "elaborate".  \\
According to ([21, 24]), the respective dimensions of the Janet bundles is known at once from the last Janet tabular by counting the number $3 + (2 \times 6) + (3 \times 6 )= 33$ of single $\bullet$, the number $6 + (3 \times 6) = 24$ of possible double $\bullet \bullet$ and finally the number $6$ of triple $\bullet \bullet \bullet$. 
  \[  \begin{array}{rccccccccccccccl}
 &&&&& 0 &&0&&0&  & 0 & & \\
 &&&&& \downarrow && \downarrow && \downarrow  & & \downarrow  &  &\\
  & 0& \longrightarrow& \Theta &\stackrel{j_3}{\longrightarrow}& 4 &\stackrel{D_1}{\longrightarrow}& 12 &\stackrel{D_2}{\longrightarrow} &  12 &\stackrel{D_3}{\longrightarrow}  & 4 & \longrightarrow  & 0  \\
  &&&&& \downarrow & & \downarrow & & \downarrow & & \downarrow &   & \\
   & 0 & \longrightarrow &  1 & \stackrel{j_3}{\longrightarrow} & 20  & \stackrel{D_1}{\longrightarrow} & 
   45  &\stackrel{D_2}{\longrightarrow} & 36 &  \stackrel{D_3}{\longrightarrow} & 10 & \longrightarrow  & 0 \\
   & & & \parallel && \hspace{5mm}\downarrow {\Phi}_0 & &\hspace{5mm} \downarrow {\Phi}_1 & & \hspace{5mm} \downarrow  {\Phi}_2&  & \hspace{5mm} \downarrow {\Phi}_3 & \\
   0 \longrightarrow & \Theta &\longrightarrow &  1 & \stackrel{\cal{D}}{\longrightarrow} &  16  & \stackrel{{\cal{D}}_1}{\longrightarrow} &  33  & \stackrel{{\cal{D}}_2}{\longrightarrow}&  24 & \stackrel{{\cal{D}}_3}{\longrightarrow} &   6  & \longrightarrow &   0    \\
   &&&&& \downarrow & & \downarrow & & \downarrow &  & \downarrow &    &   \\
   &&&&& 0 && 0 &&  0 &  & 0   &  &
   \end{array}     \]
As for the corresponding resolution of the differential module involved, it becomes:  
\[  0 \rightarrow D^6 \rightarrow D^{24} \rightarrow D^{33} \rightarrow D^{16} \rightarrow D \stackrel{p}{\rightarrow} M \rightarrow 0  \]
with Euler-Poincar\'{e} characteristic $rk_D(M) = 1 - 16 + 33 - 24 + 6 = 0$ in a coherent way.  \\

It is finally important to notice that, in the fundamental diagram I, $R^{(1)}_3$ {\it cannot} be replaced by $ {\pi}^3_2(R^{(1)}_3)= R^{(2)}_2$ after replacing $j_3$ by $j_2$ even though $dim (R^{(2)}_2)= dim (R^{(1)}_3)= 4$:

 \[  \left\{  \begin{array}{lcl}
y_{33}  & = & 0  \\
y_{23} + y_{12} & = & 0 \\
y_{22} & = & 0 \\
y_{13} & = & 0  \\
y_{11} & = & 0 \\
y_3  + y_1 & = & 0  
\end{array} 
\fbox{ $  \begin{array}{ccc}
1 & 2 & 3 \\
1 & 2 & \bullet  \\
1 & 2 & \bullet  \\
1 & \bullet & \bullet \\
1 & \bullet & \bullet \\
\bullet & \bullet & \bullet 
\end{array} $  }  \right.  \]
because $y_{112} = 0$ is missing when using the Janet tabular, the true reason for which the definition of Pommaret bases existing today in the literature is far from being intrinsic as we said in the Introduction (See [37], Example 3.14, p 119 for a similar situation showing out the importance of Spencer $\delta$-acyclicity). In the present situation, as we already said, the Spencer sequence is nothing else than the Janet sequence for the involutive first order system $R^{(1)}_3 \subset J_1(R^{(2)}_2)$. \\

In order to help the reader not familiar with the Spencer sequence, let us set $R'_1 = R^{(1)}_3$, $E' = R^{(2)}_2$ and exhibit the Janet sequence for the first order system 
$R'_1 \subset J_1(E')$ with isomorphic differential modules $M' \simeq M$ as follows (See [21], Chapter 3, Section 3 and [34] for details). The first step is to replace the four parametric jets $(y, y_1, y_2, y_{12})$ of $R^{(2)}_2$ by the new four unknowns $(z^1, z^2, z^3, z^4)$ and to consider the new first order finite type linear system $d_i z + A_i z=0$ in matrix form, obtained by substitution while taking into account $y_{112} = z^4_1= - z^4_3=0$, namely:  \\
\[  \left\{  \begin{array}{lcl}
z^1_3 + z^2 & = & 0 \\
z^2_3 & = & 0  \\
z^3_3 + z^4 & = & 0  \\
z^4_3 & = & 0  \\
z^1_2 - z^3 & = & 0  \\
z^2_2 - z^4 & = & 0  \\
z^3_2 & = & 0  \\
z^4_2 & = & 0  \\
z^1_1 - z^2 & = & 0  \\
z^2_1 & = & 0  \\
z^3_1 - z^4 & = & 0  \\
z^4_1 & = & 0  \\
\end{array}
\fbox{  $  \begin{array}{ccc}
1 & 2 & 3  \\
1 & 2 & 3  \\
1 & 2 & 3  \\
1 & 2 & 3  \\
1 & 2 & \bullet  \\
1 & 2 & \bullet  \\
1 & 2 & \bullet  \\
1 & 2 & \bullet  \\
1 & \bullet & \bullet  \\
1 & \bullet & \bullet  \\
1 & \bullet & \bullet  \\
1 & \bullet & \bullet  
\end{array}  $ } \right.  \]
The Janet tabular for this "{\it lowered system}" with $12$ unknowns has $4 + (2 \times 4) = 12$ single $\bullet$ and $4$ double $\bullet$ in a coherent way as follows:  
\[ 0 \longrightarrow  {\Theta}' \longrightarrow  E' \underset 1{\stackrel{{\cal{D}}'}{\longrightarrow}} F'_0 \underset 1{\stackrel{{\cal{D}}'_1}{\longrightarrow}} F'_1 \underset 1{\stackrel{{\cal{D}}'_2}{\longrightarrow}} F'_2 \longrightarrow 0   \]
\[  0 \longrightarrow {\Theta}' \longrightarrow 4 \longrightarrow 12 \longrightarrow 12 \longrightarrow 4 \longrightarrow 0  \]
The corresponding resolution of the differential module $M'$ is thus:
\[  0 \rightarrow D^4 \rightarrow D^{12} \rightarrow D^{12} \rightarrow D^4 \stackrel{p}{\rightarrow } M' \rightarrow \]
with Euler-Poincar\'{e} characteristic $rk_D(M') = 4 - 12 + 12 - 4 = 0$ in a coherent way.  \\ 

\noindent
{\bf EXAMPLE  2.20}: (We only provide a few hints) With the notations of the previous examples, let us consider the ideal $\mathfrak{a}=(x^3 y, y^3) \subset \mathbb{Q}[x, y] = A $. We have $rad(\mathfrak{a}) = (y) = \mathfrak{p}$ as a prime ideal. As before, we may associate the fourth order system $ R_4 \subset J_4(E)$ defined by the two PD equations $(y_{1112}=0, y_{222}=0) $ which is neither involutive nor even formally integrable. However, its symbol $g_4$, which is simply defined by the single equation $y_{1112}=0$ is trivially involutive and we have therefore ${\rho}_r(R^{(1)}_4)= R^{(1)}_{r+4}, \forall r\geq 0$. The involutive system $R^{(1)}_5$ is:
\[    \left\{  \begin{array}{lcl}
y_{22222}  & =   &  0  \\
y_{12222}  & =  &  0 \\
y_{11222} & = & 0  \\
y_{11122} & = & 0  \\
y_{11112}  &  = & 0 \\
y_{2222}  & = & 0  \\
y_{1222} & = &  0 \\
y_{1112} & = & 0 \\
y_{222} & = & 0 
\end{array} 
\fbox{ $  \begin{array}{cc}
1 & 2 \\
1 & \bullet  \\
1 & \bullet  \\
1 & \bullet \\
1 & \bullet  \\
 \bullet & \bullet  \\
 \bullet & \bullet  \\
 \bullet & \bullet \\
 \bullet & \bullet 
\end{array} $ }   \right. \]
We obtain the fundamental diagram I:
 \[  \begin{array}{rccccccccccccl}
 &&&&& 0 &&0&&0&  &  \\
 &&&&& \downarrow && \downarrow && \downarrow  & &  \\
  & 0& \longrightarrow& \Theta &\stackrel{j_3}{\longrightarrow}& C_0  &\stackrel{D_1}{\longrightarrow}&  C_1 &\stackrel{D_2}{\longrightarrow} &  C_2 & \longrightarrow  & 0  \\
  &&&&& \downarrow & & \downarrow & & \downarrow &  & \\
   & 0 & \longrightarrow &  E & \stackrel{j_3}{\longrightarrow} & C_0(E)  & \stackrel{D_1}{\longrightarrow} &  C_1(E)  &\stackrel{D_2}{\longrightarrow} & C_2(E) & \longrightarrow  & 0 \\
   & & & \parallel && \hspace{5mm}\downarrow {\Phi}_0 & &\hspace{5mm} \downarrow {\Phi}_1 & & \hspace{5mm}\downarrow {\Phi}_2 &  & \\
   0 \longrightarrow & \Theta &\longrightarrow &  E & \stackrel{\cal{D}}{\longrightarrow} &  F_0  & \stackrel{{\cal{D}}_1}{\longrightarrow} &  F_1  & \stackrel{{\cal{D}}_2}{\longrightarrow} &  F_2 & \longrightarrow  &  0      \\
   &&&&& \downarrow & & \downarrow & & \downarrow &    &   \\
   &&&&& 0 && 0 && 0  &  &
   \end{array}     \]

  \[  \begin{array}{rccccccccccccl}
 &&&&& 0 &&0&&0&  & \\
 &&&&& \downarrow && \downarrow && \downarrow  &  &\\
  & 0& \longrightarrow& \Theta &\stackrel{j_5}{\longrightarrow}& 12 &\stackrel{D_1}{\longrightarrow}& 23 &\stackrel{D_2}{\longrightarrow} &  11 &  \longrightarrow  & 0  \\
  &&&&& \downarrow & & \downarrow & & \downarrow &   & \\
   & 0 & \longrightarrow &  1 & \stackrel{j_5}{\longrightarrow} & 21  & \stackrel{D_1}{\longrightarrow} & 35  &\stackrel{D_2}{\longrightarrow} & 15 &  \longrightarrow  & 0 \\
   & & & \parallel && \hspace{5mm}\downarrow {\Phi}_0 & &\hspace{5mm} \downarrow {\Phi}_1 & & \hspace{5mm} \downarrow  {\Phi}_2&  &  \\
   0 \longrightarrow & \Theta &\longrightarrow &  1 & \stackrel{\cal{D}}{\longrightarrow} &  9  & \stackrel{{\cal{D}}_1}{\longrightarrow} &  12  & \stackrel{{\cal{D}}_2}{\longrightarrow}&  4 &  \longrightarrow &   0    \\
   &&&&& \downarrow & & \downarrow & & \downarrow &  &  \\
   &&&&& 0 && 0 &&  0 &  & 
   \end{array}     \]
We have in particular the short exact sequence:
\[  0 \rightarrow J_6(E) \rightarrow J_1(J_5(E)) \rightarrow C_1(E) \rightarrow   0  \,\,\,  \Rightarrow  \,\,\, 0 \rightarrow  28  \rightarrow 63 \rightarrow  35 \rightarrow 0  \]
The comparison between the use of the Koszul homology and the use of the Spencer sequence needs again no comment but the reader will appreciate trying to do  it in an intrinsic way by making any linear change of the coordinates.  \\
The final hint is to introduce the symbol $g'_5 \subset S_5 T^* \otimes E$ of $R'_5 = R^{(1)}_5$ and its prolongation $g'_6 \subset S_6 T^* \otimes E$ in the definition of the Spencer bundles with now $C_0 = R'_5, C_1 = T^* \otimes R'_5 / \delta g'_6, C_2 = {\wedge}^2 T^* \otimes R'_5 / \delta (T^* \otimes g'_6)   $, getting for example 
$dim(C_1) = (2 \times 12) - 1 = 23 $ and $dim(C_2)= 12 - 1 = 11$ through the short exact $\delta$-sequences like: 
\[ 0 \rightarrow g'_7 \stackrel{\delta}{\rightarrow} T^*\otimes g'_6 \stackrel{\delta}{\rightarrow} {\wedge}^2 T^* \otimes g'_5 \rightarrow 0 \,\, \Rightarrow \,\, 0 \rightarrow 1 \rightarrow  2  \rightarrow  1 \rightarrow 0  \]  
Once again, no classical method can provide these results.  \\

\noindent
{\bf EXAMPLE 2.21}: We study in an intrinsic way an example proposed by V. Gerdt in $2000$ ([9]). Like in the previous examples, let us consider the ideal $\mathfrak{a} = (x^2 y - z, x y^2 - y) \subset \mathbb{Q} [x, y, z]$. We transform it into the third order system defined by the two corresponding PD equations, namely $(y_{112} - y_3=0, y_{122} - y_2 = 0)$. This system is neither involutive, nor even formally integrable. Using crossed derivatives, we obtain two new second order PD equations 
$(y_{23} - y_{12}=0, y_{33} - y_{13}=0)$ because $ y_{33} = y_{1123} = y_{1112} = y_{13}$ after two prolongations that we may differentiate once more in order to get the two new third order PD equations $ (y_{223} - y_2 =0, y_{233} - y_3 = 0)$. Accordingly, as $y_{233}= y_{123}$, we have thus obtained an equivalent third order system $R_3 \subset J_3(E)$ defined by the following PD equations:
 \[  \left\{  \begin{array}{lcl}
y_{333} - y_{113}  & =   &  0  \\
y_{233} - y_3  & =  &  0 \\
y_{223} - y_2 & = & 0  \\
y_{133} - y_{113}  &  = & 0 \\
y_{123} - y_3  & = & 0  \\
y_{122} - y_2 & = &  0 \\
y_{112} - y_3  & = & 0 \\
y_{33} - y_{13}  & = & 0  \\
y_{23} - y_{12} & = & 0 
\end{array} 
\fbox{ $  \begin{array}{ccc}
1 & 2 & 3 \\
1 & 2 & \bullet  \\
1 & 2 & \bullet  \\
1 & \bullet & \bullet  \\
1 & \bullet & \bullet  \\
1 & \bullet & \bullet  \\
1 & \bullet & \bullet \\
\bullet & \bullet & \bullet \\
\bullet & \bullet & \bullet 
\end{array} $ }   \right. \]
This system may not be involutive because $d_2 (y_{122} - y_2) $ cannot vanish but we may modify the ordering by changing coordinates. As the character ${\alpha}^2_3$ must be zero, we must in any case make the class $2$ full by introducing $y_{222}$ as leading term. The easiest possibility is to make the change of coordinates 
$({\bar{x}}^1 = x^1, {\bar{x}}^2 = x^2 + x^1 , {\bar{x}}^3 = x^3)$. Suppressing the bar, we obtain the new equations:  \\ 
\[  \left\{  \begin{array}{lcl}
y_{333} - y_{113} + y_2 - 2 y_3  & =   &  0  \\
y_{233} - y_3  & =  &  0 \\
y_{223} - y_2 & = & 0  \\
y_{222} - y_{112}  + y_3 - 2 y_2  & = & 0  \\
y_{133} - y_{113} + y_2 - y_3  &  = & 0  \\
y_{123} + y_2 - y_3  & = & 0  \\
y_{122}  + y_{112} + y_2- y_3 & = &  0 \\
y_{33} - y_{23} - y_{13}   & = & 0  \\
y_{23} - y_{22} -  y_{12}  & = & 0 
\end{array} 
\fbox{ $  \begin{array}{ccc}
1 & 2 & 3 \\
1 & 2 & \bullet  \\
1 & 2 & \bullet  \\
1 & 2 & \bullet  \\
1 & \bullet & \bullet  \\
1 & \bullet & \bullet  \\
1 & \bullet & \bullet \\
\bullet & \bullet & \bullet \\
\bullet & \bullet & \bullet 
\end{array} $ }   \right. \]
The system has an involutive symbol if we set $y_2 = 0, y_3 = 0$ and is formally integrable, thus involutive in this new coordinates. {\it As it is an intrinsic property, it is thus involutive in any coordinate system}. It follows that we have at once the Janet sequence by counting $15$ single $\bullet$, $3 + (2 \times 3) = 9$ double $\bullet \bullet$ and $2$ triple $\bullet \bullet \bullet$. We have thus $dim(E) = 1, dim(F_0) = 9, dim(F_1) = 15, dim(F_2) = 9, dim(F_3) = 2$. The fundamental diagram I becomes:

  \[  \begin{array}{rccccccccccccccl}
 &&&&& 0 &&0&&0&  & 0 & & \\
 &&&&& \downarrow && \downarrow && \downarrow  & & \downarrow  &  &\\
  & 0& \longrightarrow& \Theta &\stackrel{j_3}{\longrightarrow}& 11 &\stackrel{D_1}{\longrightarrow}& 30 &\stackrel{D_2}{\longrightarrow} &  27 &\stackrel{D_3}{\longrightarrow}  & 8 & \longrightarrow  & 0  \\
  &&&&& \downarrow & & \downarrow & & \downarrow & & \downarrow &   & \\
   & 0 & \longrightarrow &  1 & \stackrel{j_3}{\longrightarrow} & 20  & \stackrel{D_1}{\longrightarrow} & 
   45  &\stackrel{D_2}{\longrightarrow} & 36 &  \stackrel{D_3}{\longrightarrow} & 10 & \longrightarrow  & 0 \\
   & & & \parallel && \hspace{5mm}\downarrow {\Phi}_0 & &\hspace{5mm} \downarrow {\Phi}_1 & & \hspace{5mm} \downarrow  {\Phi}_2&  & \hspace{5mm} \downarrow {\Phi}_3 & \\
   0 \longrightarrow & \Theta &\longrightarrow &  1 & \stackrel{\cal{D}}{\longrightarrow} &  9 & \stackrel{{\cal{D}}_1}{\longrightarrow} &  15  & \stackrel{{\cal{D}}_2}{\longrightarrow}&  9 & \stackrel{{\cal{D}}_3}{\longrightarrow} &   2  & \longrightarrow &   0    \\
   &&&&& \downarrow & & \downarrow & & \downarrow &  & \downarrow &    &   \\
   &&&&& 0 && 0 &&  0 &  & 0   &  &
   \end{array}     \]
It is {\it much more difficult} to check the dimensions of the Spencer bundles. Calling again $R_3 \subset J_3(E)$ this involutive system, the characters of its symbol $g_3$ are ${\alpha}^3_3 = 0, {\alpha}^2_3 = 0, {\alpha}^1_3 = 6 - 3 = 3$. We obtain $dim(g_4) = (1 \times 3) + (2 \times 0) + (3 \times 0) =  3$ and $dim(g_{r+3})=3, \forall r\geq 0$ in the exact $\delta$-sequence:
\[      0 \rightarrow g_5 \stackrel{\delta}{\rightarrow} T^* \otimes g_4 \stackrel{\delta}{\rightarrow} {\wedge}^2 T^* \otimes g_3 \stackrel{\delta}{ \rightarrow} {\wedge}^3 T^* \otimes S_2 T^* \otimes E  \]  
we obtain for example $C_2 : {\wedge}^2 T^* \otimes R_3/ \delta (T^* \otimes g_4) \Rightarrow dim(C_2) = (3 \times 11) - ((3\times 3) - 3) = 33 - 6 = 27$ in a coherent way. Once more, {\it no classical method could provide these results}. \\

\noindent
{\bf EXAMPLE  2.22}: ({\it Macaulay} in [18]) Among the best examples we know that justify our comments on Pommaret bases, we shall revisit one dealing with formal integrability and one, on the contrary, which is trivially FI but is among the few rare elementary explicit examples of a $2$-acyclic symbol which is not involutive, apart from the symbol of the conformal Killing operator for a non-degenerate metric that we shall consider later on. \\
\noindent
1) With the notations of the previous examples, let us consider the ideal $\mathfrak{a}= (z^3, xz - y) \subset k[x, y, z]$ with $rad(\mathfrak{a})= (y, z) = \mathfrak{p}$ a prime ideal. We may transform it into a non-homogeneous second order system of PD equations $R_2 \subset J_2(E)$ defined by $(y_{33} = 0, y_{13} - y_2 = 0)$. Exchanging $1$ with $2$, the symbol $g_2$ defined by $(y_{33} = 0, y_{23} = 0)$ is involutive. We have therefore ${\rho}_r(R^{(1)}_2 = R^{(1)}_{r+2}$ and in particular $ R^{(1)}_2$ is defined by $(y_{33}=0, y_{23}=0, y_{13} - y_2=0)$. By chance, its symbol $g^{(2)}_2$ defined by $(y_{33}=0, y_{23}=0, y_{13}=0)$ is {\it again} involutive with one equation of class $1$, one equation of class $2$ and one equation of class $1$. Accordingly, ${\rho}_1(R^{(1)}_2)=R^{(1)}_3$ projects onto $R^{(2)}_2$ defined by $(y_{33}=0, y_{23}=0, y_{22}=0, y_{13} - y_2=0)$ and the PP procedure is ending because this system is involutive: \\
 \[  \left\{  \begin{array}{lcl}
y_{33}   & =   &  0  \\
y_{23}  & =  &  0 \\
y_{22}  & = & 0  \\
y_{13} - y_2  &  = & 0 
\end{array} 
\fbox{ $  \begin{array}{ccc}
1 & 2 & 3 \\
1 & 2 & \bullet  \\
1 & 2 & \bullet  \\
1 & \bullet & \bullet  
\end{array} $ }   \right. \]
We have the strict inclusions $ R'_2 = R^{(2)}_2 \subset R^{(1)}_2 \subset R_2 \subset J_2(E)$ and the fundamental diagram I:  \\
\[  \begin{array}{rccccccccccccccl}
 &&&&& 0 &&0&&0&  & 0 & & \\
 &&&&& \downarrow && \downarrow && \downarrow  & & \downarrow  &  &\\
  & 0& \longrightarrow& \Theta &\stackrel{j_2}{\longrightarrow}& 6 &\stackrel{D_1}{\longrightarrow}& 16 &\stackrel{D_2}{\longrightarrow} &  14 &\stackrel{D_3}{\longrightarrow}  & 4 & \longrightarrow  & 0  \\
  &&&&& \downarrow & & \downarrow & & \downarrow & & \downarrow &   & \\
   & 0 & \longrightarrow &  1 & \stackrel{j_2}{\longrightarrow} & 10  & \stackrel{D_1}{\longrightarrow} & 
   20  &\stackrel{D_2}{\longrightarrow} & 15 &  \stackrel{D_3}{\longrightarrow} & 4 & \longrightarrow  & 0 \\
   & & & \parallel && \hspace{5mm}\downarrow {\Phi}_0 & &\hspace{5mm} \downarrow {\Phi}_1 & &  \downarrow  &  &  \downarrow  & \\
   0 \longrightarrow & \Theta &\longrightarrow &  1 & \stackrel{\cal{D}}{\longrightarrow} &  4  & \stackrel{{\cal{D}}_1}{\longrightarrow} &  4  & \stackrel{{\cal{D}}_2}{\longrightarrow} &  1
   & \longrightarrow  &  0  & &       \\
   &&&&& \downarrow & & \downarrow & & \downarrow  &  & &    &   \\
   &&&&& 0 && 0 && 0  &  &    &  &
   \end{array}     \]
   In fact, we have the final characters ${\alpha}^3_2 = 0, {\alpha}^2_2 = 0, {\alpha}^1_2=3-1=2$ and thus successively:
   \[ C'_0= R'_2 , C'_1 = T^* \otimes R'_2 / \delta (g'_3), C'_2 = {\wedge}^2 T^* \otimes R'_2 / \delta (T^* \otimes g'_3), 
   C'_3 = {\wedge}^3 T^* \otimes R'_2 / \delta ({\wedge}^2 T^* \otimes g'_3) \]
with $dim(R'_2) = 6$ and $dim(g'_{r+2}) = 2, \forall r \geq 0$ in the exact $\delta$-sequence:  \\
\[ 0 \rightarrow g'_4 \stackrel{\delta}{\rightarrow } T^* \otimes g'_3 \stackrel{\delta}{\rightarrow} {\wedge}^2 T^* \otimes g'_2 \stackrel{\delta}{\rightarrow} {\wedge}^3 T^* \otimes T^*  \]
We check indeed:  \\
 \[ dim(C'_0)=6, dim(C'_1)=(3\times 6) - 2=16, dim(C'_2)=(3 \times 6) - (( 3 \times 2) - 2) = 18 - 4 = 14, \]  
 \[  dim(C'_3)= 6 - dim({\wedge}^3 T^*  \otimes g'_2)=6-2=4 \]
Again, no classical method could provide such results.

\noindent
2) With the notations of the previous examples, let us consider the homogeneous ideal $\mathfrak{a}= (z^2, y z - x^2, y^2) \subset k[x, y, z]$ with $rad(\mathfrak{a})= (x, y, z)=\mathfrak{m}$ a maximal and thus zero-dimensional prime ideal. We may transform it into an homogeneous second order system of PD equations $R_2 \subset J_2(E)$ defined by 
$(y_{33}= 0, y_{23} - y_{11}=0, y_{22}= 0)$ but the reader may treat as well the system $(y_{33} - y_{11}=0, y_{23}=0, y_{22} - y_{11}=0)$. Of course, this system is FI because it is homogeneous but we let the reader check on the Janet tabular that $g_2$ is not involutive though the coordinate system is surely $\delta$-regular because be have full class $3$ and full class $2$. All the third order jets vanish but $y_{123} - y_{111}=0$ leading to $dim (g_3)=1 \Rightarrow dim (R_3) = 8 = 2^3$ ([18]). Finally $g'_4=0 \Rightarrow dim(R_4)= 8$ and we could believe that we do not need any PP procedure as $R_4$ is an ivolutive system because $g_4$ is trivially involutive and $R_2$ is finite type like the killing system. Moreover, as we have constant coefficients, the three brackets and their only Jacobi identity do provide the following sequence which is of course quite far from being a Janet sequence as it only involves second order operators. Also, {\it It is important to notice that the knowledge of the first second order operator does not provide any way to obtain the third without passing through the second, contrary to the situation existing in the Janet sequence}:  \\ 
\[  \begin{array}{rcccccccccccl}
0 & \rightarrow & \Theta & \rightarrow & E & \underset 2{\rightarrow} & F_0 & \underset 2{\rightarrow} & F_1 & \underset 2{\rightarrow}& F_2 & \rightarrow   &  0  \\

0 & \rightarrow & \Theta & \rightarrow & 1 & \rightarrow & 3 & \rightarrow & 3 & \rightarrow & 1 & \rightarrow  & 0  
\end{array}  \]
Moreover, such a procedure is rather "experimental" and must be coherent with the theorem saying that the order of generating CC is one plus the number of prolongations needed to reach a $2$-acyclic symbol, that is $g_3$ {\it must be} $2$-acyclic. Equivalently the $\delta$-sequence: \\
\[   0 \rightarrow  {\wedge}^2 T^* \otimes g_3 \stackrel{\delta}{\rightarrow} {\wedge}^3 T^* \otimes g_2 \rightarrow 0  \] 
must be exact. We let the reader prove that the corresponding $3 \times 3$ matrix has maximum rank or refer to ([  ]) for an explicit computation.  \\
It remains to work out the corresponding Janet and Spencer sequences and there is a first delicate point to overcome. Indeed, as $g_4=0$, we have thus $R_4 \simeq R_3$ and it could be tempting to start with $R_3$ and thus to replace $j_4$ by $j_3$ in the fundamental diagram I. However, it should lead to a dead end because $C_0 = R_q \subset J_q(E)$ {\it must be} an involutive system (See [37], Example  3.14, p 119 to 126 ). Accordingly, using $j_4$ and the fact that $g_4=0 \Rightarrow g_5=0$, we have $C_r = {\wedge}^r T^* \otimes R_4, \forall r=0, 1, 2, 3$.The interest of this approach, having no meaning for anybody not aware of the Spencer operator $d$ and its symbol restriction $\delta$, is that we can obtain the Janet sequence without any explicit calculation as follows through the formula $dim(C_r(E))= dim (C_r) + dim (F_r)$ in the following diagram:  \\

  \[  \begin{array}{rccccccccccccccl}
 &&&&& 0 &&0&&0&  & 0 & & \\
 &&&&& \downarrow && \downarrow && \downarrow  & & \downarrow  &  &\\
  & 0& \longrightarrow& \Theta &\stackrel{j_4}{\longrightarrow}& 8 &\stackrel{D_1}{\longrightarrow}& 24 &\stackrel{D_2}{\longrightarrow} &  24 &\stackrel{D_3}{\longrightarrow}  & 8 & \longrightarrow  & 0  \\
  &&&&& \downarrow & & \downarrow & & \downarrow & & \downarrow &   & \\
   & 0 & \longrightarrow &  1 & \stackrel{j_4}{\longrightarrow} & 35  & \stackrel{D_1}{\longrightarrow} & 
   84  &\stackrel{D_2}{\longrightarrow} & 70&  \stackrel{D_3}{\longrightarrow} & 20 & \longrightarrow  & 0 \\
   & & & \parallel && \hspace{5mm}\downarrow {\Phi}_0 & &\hspace{5mm} \downarrow {\Phi}_1 & & \hspace{5mm} \downarrow  {\Phi}_2&  & \hspace{5mm} \downarrow {\Phi}_3 & \\
   0 \longrightarrow & \Theta &\longrightarrow &  1 & \underset 4{ \stackrel{\cal{D}}{\longrightarrow}} &  27  & \stackrel{{\cal{D}}_1}{\longrightarrow} &  60  & \stackrel{{\cal{D}}_2}{\longrightarrow}&  46 & \stackrel{{\cal{D}}_3}{\longrightarrow} &   12  & \longrightarrow &   0    \\
   &&&&& \downarrow & & \downarrow & & \downarrow &  & \downarrow &    &   \\
   &&&&& 0 && 0 &&  0 &  & 0   &  &
   \end{array}     \]   
   
\noindent
{\bf 3) DIFFERENTIAL DUALITY}  \\

For example, the fact that the Cauchy operator is the adjoint of the Killing operator for the Euclidean metric is in any textbook of continuum mechanics in the chapter "variational calculus" and the parametrization problem has been quoted by many famous authors, as we said in the Abstract, but only from a computational point of view. It is still not known that the adjoint of the $20$ components of the Bianchi operator has been introduced by C. Lanczos as we explained with details in ([35, 40, 48]). The main trouble is that these two problems have {\it never} ben treated in an intrinsic way and, in particular, changes of coordinates have {\it never} been considered. The same situation can be met for Maxwell equations but is out of our scope [23, 33, 35, 41, 44, 47, 49].   \\

\noindent
{\bf LEMMA 3.1}: When $y^k=f^k(x)$ is invertible with $\Delta(x)=det({\partial}_i f^k(x))\neq 0$ and inverse $x=g(y)$, then we have $n$ {\it identities} $ \frac{\partial}{\partial y^k} (\frac{1}{\Delta} 
{\partial }_if ^k(g(y))) = 0 $.  \\

\noindent
{\bf PROPOSITION 3.2}: The Cauchy operator is the adjoint of the Killing operator.  \\

{\it Proof}: Let $X$ be a manifold of dimension $n$ with local coordinates $(x^1, ... , x^n)$, tangent bundle $T$ and cotangent bundle $T^*$. If $\omega \in S_2 T^*$ is a metric with $det(\omega)\neq 0$, we my introduce the standard Lie derivative in order to define the first order Killing operator:    \\
\[       {\cal{D}}: \xi \in T \rightarrow \Omega = 
({\Omega}_{ij}= {\omega}_{rj}(x) {\partial}_i{\xi}^r + {\omega}_{ir}(x){\partial}_j{\xi}^r + 
{\xi}^r {\partial}_r {\omega}_{ij}(x) ) \in S_2T^*  \]
Here start the problems because, in our opinion at least, a systematic use of the adjoint operator has never been used in mathematical physics and even in continuum mechanics apart through a variational procedure. As will be seen later on, the purely intrinsic definition of the adjoint can only be done in the theory of differential modules by means of the so-called {\it side changing functor}. From a purely differential geometric point of view, the idea is to associate to any vector bundle $E$ over $X$ a new vector bundle $ad(E)= {\wedge}^nT^* \otimes E^*$ where $E^*$ is obtained from $E$ by patching local coordinates while inverting the transition matrices, exactly like $T^*$ is obtained from $T$. It follows that the stress tensor $\sigma =({\sigma}^{ij}) \in ad(S_2T^*) =  {\wedge}^n T^* \otimes S_2 T$ is {\it not} a tensor but a tensor density, that is transforms like a tensor up to a certain power of the Jacobian matrix. When $n=4$,  the fact that such an object is called stress-energy tensor does not change anything as it cannot be related to the Einstein tensor which is a true {\it tensor} indeed. Of course, it is always possible in GR to use $(det(\omega))^\frac{1}{2} $ but, as we shall see, the study of contact structures {\it must} be done without any reference to a background metric. In any case, we may define as usual:  \\
\[  ad({\cal{D}}): {\wedge}^n T^* \otimes S_2T \rightarrow {\wedge}^n T^* \otimes T : \sigma \rightarrow \varphi \]
Multiplying ${\Omega}_{ij}$ by ${\sigma}^{ij}$ and integrating by parts, the factor of 
$ - 2 \, {\omega}_{kr}\,  {\xi}^r$ is easly seen to be:  \\
\[   {\nabla}_i {\sigma}^{ik} = {\partial}_i {\sigma}^{ik} +  {\gamma}^k_{ij} {\sigma}^{ij} = {\varphi}^k   \] 
with well known Christoffel symbols 
$ {\gamma}^k_{ij} = \frac{1}{2} {\omega}^{kr} ({\partial}_i {\omega}_{rj} + {\partial}_j {\omega}_{ir} - 
{\partial}_r {\omega}_{ij}) $. \\
However, if the stress should be a tensor, we should get for the covariant derivative:  \\
\[ {\nabla}_r {\sigma}^{ij}= {\partial}_r{\sigma}^{ij} + {\gamma}^i_{rs}{\sigma}^{sj}+ {\gamma}^j_{rs} {\sigma}^{is} \Rightarrow 
{\nabla}_i {\sigma}^{ik} = {\partial}_i {\sigma}^{ik} + {\gamma}^r_{ri} {\sigma}^{ik} + 
{\gamma}^k_{ij} {\sigma}^{ij}    \]
The difficulty is to prove that we do not have a contradiction because $\sigma$ is a tensor density.  \\
If we have an invertible transformation like in the lemma, we have successively by using it:  \\
\[ {\tau}^{kl}(f(x))= \frac{1}{\Delta} {\partial}_i f^k(x) {\partial}_j f^l(x) {\sigma}^{ij}(x) \]
\[  \frac{\partial {\tau}^{kl}}{\partial y^k}= \frac{\partial}{\partial y^k} ( (\frac{1}{\Delta} {\partial}_if^k) + 
\frac{1}{\Delta}{\partial}_i f^k  \frac{\partial}{\partial y^k}({\partial}_j f^l) {\sigma}^{ij}+ 
\frac{1}{\Delta}{\partial}_i f^k {\partial}_j f^l \frac{\partial}{\partial y^k} {\sigma}^{ij}  \]
\[  \frac{\partial {\tau}^{ku}}{\partial y^k} = \frac{1}{\Delta} ({\partial}_{ij} f^u) {\sigma}^{ij}  + \frac{1}{\Delta} {\partial}_j f^u {\partial}_i {\sigma}^{ij}   \]
Now, we recall the transformation law of the Christoffel symbols, namely:   \\
\[    {\partial}_r f^u(x) {\gamma}^r_{ij} (x) =  {\partial}_{ij} f^u(x)  + {\partial}_i f^k(x) {\partial}_j f^l(x)  {\bar{\gamma}}^u_{kl} (f(x))  \]
\[ \Rightarrow   \frac{1}{\Delta}{\partial}_r f^u {\gamma}^r_{ij} {\sigma}^{ij}= \frac{1}{\Delta} {\partial}_{ij} f^u {\sigma}^{ij} + {\bar{\gamma}}^u_{kl}(y) {\tau}^{kl}   \]
Eliminating the second derivatives of $f$ we finally get:  \\
\[   {\psi}^u= \frac{\partial {\tau}^{ku}}{\partial y^k} + {\bar{ \gamma}}^u_{kl} = \frac{1}{\Delta} {\partial}_r f^u({\partial}_i {\sigma}^{ir} + {\gamma}^r_{ij} {\sigma}^{ij}) = \frac{1}{\Delta} {\partial}_r f^u {\varphi}^r  \]      
This tricky technical result, which is not evident at all, explains why the additional term we had is just disappearing in fact when $\sigma$ is a density. \\
\hspace*{12cm}  $   \Box   $   \\

Let $K$ be a {\it differential field} with $n$ commuting derivations $({\partial}_1,...,{\partial}_n)$ and consider the ring $D=K[d_1,...,d_n]=K[d]$ of differential operators with coefficients in $K$ with $n$ commuting formal derivatives satisfying $d_ia=ad_i + {\partial}_ia$ in the operator sense. If $P=a^{\mu}d_{\mu}\in D=K[d]$, the highest value of ${\mid}\mu {\mid}$ with $a^{\mu}\neq 0$ is called the {\it order} of the {\it operator} $P$ and the ring $D$ with multiplication $(P,Q)\longrightarrow P\circ Q=PQ$ is filtred by the order $q$ of the operators. We have the {\it filtration} $0\subset K=D_0\subset D_1\subset  ... \subset D_q \subset ... \subset D_{\infty}=D$. As an algebra, $D$ is generated by $K=D_0$ and $T=D_1/D_0$ with $D_1=K\oplus T$ if we identify an element $\xi={\xi}^id_i\in T$ with the vector field $\xi={\xi}^i(x){\partial}_i$ of differential geometry, but with ${\xi}^i\in K$ now. It follows that $D={ }_DD_D$ is a {\it bimodule} over itself, being at the same time a left $D$-module by the composition $P \longrightarrow QP$ and a right $D$-module by the composition $P \longrightarrow PQ$. We define the {\it adjoint} functor $ad:D \longrightarrow D^{op}:P=a^{\mu}d_{\mu} \longrightarrow  ad(P)=(-1)^{\mid \mu \mid}d_{\mu}a^{\mu}$ and we have $ad(ad(P))=P$ both with $ad(PQ)=ad(Q)ad(P), \forall P,Q\in D$. Such a definition can be extended to any matrix of operators by using the transposed matrix of adjoint operators (See [16, 27, 28, 38, 39, 42, 53] for more details and applications to control theory or mathematical physics). \\

Accordingly, if $y=(y^1, ... ,y^m)$ are differential indeterminates, then $D$ acts on $y^k$ by setting $d_iy^k=y^k_i \longrightarrow d_{\mu}y^k=y^k_{\mu}$ with $d_iy^k_{\mu}=y^k_{\mu+1_i}$ and $y^k_0=y^k$. We may therefore use the jet coordinates in a formal way as in the previous section. Therefore, if a system of OD/PD equations is written in the form ${\Phi}^{\tau}\equiv a^{\tau\mu}_ky^k_{\mu}=0$ with coefficients $a\in K$, we may introduce the free differential module $Dy=Dy^1+ ... +Dy^m\simeq D^m$ and consider the differential {\it module of equations} $I=D\Phi\subset Dy$, both with the residual {\it differential module} $M=Dy/D\Phi$ or $D$-module and we may set $M={ }_DM$ if we want to specify the ring of differential operators. We may introduce the formal {\it prolongation} with respect to $d_i$ by setting $d_i{\Phi}^{\tau}\equiv a^{\tau\mu}_ky^k_{\mu+1_i}+({\partial}_ia^{\tau\mu}_k)y^k_{\mu}$ in order to induce maps $d_i:M \longrightarrow M:{\bar{y} }^k_{\mu} \longrightarrow {\bar{y}}^k_{\mu+1_i}$ by residue with respect to $I$ if we use to denote the residue $Dy \longrightarrow M: y^k \longrightarrow {\bar{y}}^k$ by a bar like in algebraic geometry. However, for simplicity, we shall not write down the bar when the background will indicate clearly if we are in $Dy$ or in $M$. As a byproduct, the differential modules we shall consider will always be {\it finitely generated} ($k=1,...,m<\infty$) and {\it finitely presented} ($\tau=1, ... ,p<\infty$). Equivalently, introducing the {\it matrix of operators} ${\cal{D}}=(a^{\tau\mu}_kd_{\mu})$ with $m$ columns and $p$ rows, we may introduce the morphism $D^p \stackrel{{\cal{D}}}{\longrightarrow} D^m:(P_{\tau}) \longrightarrow (P_{\tau}{\Phi}^{\tau})$ over $D$ by acting with $D$ {\it on the left of these row vectors} while acting with ${\cal{D}}$ {\it on the right of these row vectors} by composition of operators with $im({\cal{D}})=I$. The {\it presentation} of $M$ is defined by the exact cokernel sequence $D^p \stackrel{{\cal{D}}}{\longrightarrow} D^m \longrightarrow M \longrightarrow 0 $. We notice that the presentation only depends on $K, D$ and $\Phi$ or $ \cal{D}$, that is to say never refers to the concept of (explicit local or formal) solutions. It follows from its definition that $M$ can be endowed with a {\it quotient filtration} obtained from that of $D^m$ which is defined by the order of the jet coordinates $y_q$ in $D_qy$. We have therefore the {\it inductive limit} $0 \subseteq M_0 \subseteq M_1 \subseteq ... \subseteq M_q \subseteq ... \subseteq M_{\infty}=M$ with $d_iM_q\subseteq M_{q+1}$ and $M=DM_q$ for $q\gg 0$ with prolongations $D_rM_q\subseteq M_{q+r}, \forall q,r\geq 0$. It is important to notice that it may be sometimes quite difficult to work out $I_q$ or $M_q$ from a given presentation which is not involutive [24, 28].   \\ 

\noindent
{\bf DEFINITION 3.3}: An exact sequence of morphisms finishing at $M$ is said to be a {\it resolution} of $M$. If the differential modules involved apart from $M$ are free, that is isomorphic to a certain power of $D$, we shall say that we have a {\it free resolution} of $M$.  \\

Having in mind that $K$ is a left $D$-module with the action $(D,K) \longrightarrow K:(d_i,a)\longrightarrow {\partial}_ia$ and that $D$ is a bimodule over itself with $PQ\neq QP$, {\it we have only two possible constructions}:  \\

\noindent
{\bf DEFINITION 3.4}: We may define the right ({\it care})  differential module $hom_D(M,D) $ with $(f P)(m) = (f(m)) P \Rightarrow (fPQ)(m)=((fP)(m))Q=((f(m))P)Q=(f(m))PQ$.  \\                        

\noindent
{\bf DEFINITION 3.5}: We define the {\it system} $R=hom_K(M,K)$ and set $R_q=hom_K(M_q,K)$ as the {\it system of order} $q$. We have the {\it projective limit} $R=R_{\infty} \longrightarrow ... \longrightarrow R_q \longrightarrow ... \longrightarrow R_1 \longrightarrow R_0$. It follows that $f_q\in R_q:y^k_{\mu} \longrightarrow f^k_{\mu}\in K$ with $a^{\tau\mu}_kf^k_{\mu}=0$ defines a {\it section at order} $q$ and we may set $f_{\infty}=f\in R$ for a {\it section} of $R$. For an arbitrary differential field $K$, {\it such a definition has nothing to do with the concept of a formal power series solution} ({\it care}).\\

\noindent
{\bf PROPOSITION 3.6}: When $M$ is a left $D$-module, then $R$ is also a left $D$-module. \\

\noindent
{\it Proof}: As $D$ is generated by $K$ and $T$ as we already said, let us define:  \\
\[  (af)(m)=af(m)= f(am), \hspace{4mm} \forall a\in K, \forall m\in M \]
\[ (\xi f)(m)=\xi f(m)-f(\xi m), \hspace{4mm} \forall \xi=a^id_i\in T,\forall m\in M  \]
In the operator sense, it is easy to check that $d_ia=ad_i+{\partial}_ia$ and that $\xi\eta - \eta\xi=[\xi,\eta]$ is the standard bracket of vector fields. We finally 
get $(d_if)^k_{\mu}=(d_if)(y^k_{\mu})={\partial}_if^k_{\mu}-f^k_{\mu +1_i}$ and thus recover {\it exactly} the Spencer operator of the previous section though {\it this is not evident at all}. We also get $(d_id_jf)^k_{\mu}={\partial}_{ij}f^k_{\mu}-{\partial}_if^k_{\mu+1_j}-{\partial}_jf^k_{\mu+1_i}+f^k_{\mu+1_i+1_j} \Longrightarrow d_id_j=d_jd_i, \forall i,j=1,...,n$ and thus $d_iR_{q+1}\subseteq R_q\Longrightarrow d_iR\subset R$ induces a well defined operator $R\longrightarrow T^*\otimes R:f \longrightarrow dx^i\otimes d_if$. This operator has been first introduced, up to sign, by F.S. Macaulay as early as in $1916$ but this is still not ackowledged [18, 30]. (See  [31, 38, 42] for more details and applications).  \\
\hspace*{12cm}   $  \Box   $   \\

\noindent
{\bf DEFINITION 3.7}: With any differential module $M$ we shall associate the {\it graded module} $G=gr(M)$ over the polynomial ring $gr(D)\simeq K[\chi]$ by setting $G={\oplus}^{\infty}_{q=0} G_q$ with $G_q=M_q/M_{q-1}$ and we get $g_q=G_q^*$ where the {\it symbol} $g_q$ is defined by the short exact sequences: \\
\[ 0\longrightarrow M_{q-1}\longrightarrow M_q \longrightarrow G_q \longrightarrow 0  \hspace{4mm}  \Longleftrightarrow \hspace{4mm}  0 \longrightarrow g_q \longrightarrow R_q \longrightarrow R_{q-1} \longrightarrow 0  \]
We have the short exact sequences $0\longrightarrow D_{q-1} \longrightarrow D_q \longrightarrow S_qT \longrightarrow 0 $ leading to $gr_q(D)\simeq S_qT$ and we may set as usual $T^*=hom_K(T,K)$ in a coherent way with differential geometry. \\

The two following definitions, which are well known in commutative algebra, are also valid (with more work) in the case of differential modules (See [27, 28] for more details or the references [16, 20, 27, 28, 52] for an introduction to homological algebra and diagram chasing).  \\

\noindent
{\bf DEFINITION 3.8}: The set of elements $t(M) = \{ m \in M \mid \exists 0 \neq P \in D, Pm=0\}\subseteq M$ is a differential module called the {\it torsion submodule} of $M$. More generally, a module $M$ is called a {\it torsion module} if $t(M)=M$ and a {\it torsion-free module} if $t(M)=0$. In the short exact sequence $0 \rightarrow t(M) \rightarrow M \rightarrow M' \rightarrow 0$, the module $M'$ is torsion-free. Its defining module of equations $I'$ is obtained by adding to $I$ a representative basis of $t(M)$ set up to zero and we have thus $I \subseteq I'$.  \\

\noindent
{\bf DEFINITION 3.9}: A differential module $F$ is said to be {\it free} if $F \simeq D^r$ for some integer 
$r > 0$ and we shall {\it define} $rk_D(F)=r$. If $F$ is the biggest free differential module contained in $M$, then $M/F$ is a torsion differential module and $hom_D(M/F,D)=0$. In that case, we shall {\it define} the {\it differential rank} of $M$ to be $rk_D(M)=rk_D(F)=r$. Accordingly, if $M$ is defined by a linear involutive operator of order $q$, then $rk_D(M)={\alpha}^n_q$.  \\

\noindent
{\bf PROPOSITION 3.10}: If $0 \rightarrow M' \rightarrow M \rightarrow M" \rightarrow 0$ is a short exact sequence of differential modules and maps or operators, we have $rk_D(M)=rk_D(M') + rk_D(M")$.  \\

In the general situation, let us consider the sequence $ M' \stackrel{f}{\longrightarrow} M  \stackrel{g}{\longrightarrow} M" $ of modules which may not be exact and define $B =im(f) \subseteq Z = ker(g) \Rightarrow   H=Z/B$. \\

In order to conclude this section, we may say that the main difficulty met when passing from the differential framework to the algebraic framework is the " {\it inversion} " of arrows. Indeed, when an operator is injective, that is when we have the exact sequence $ 0 \rightarrow E \stackrel{{\cal{D}}}{\longrightarrow} F$ with $dim(E)=m, dim(F)=p$, like in the case of the operator $0 \rightarrow E \stackrel{j_q}{\longrightarrow} J_q(E) $, on the contrary, using differential modules, we have the epimorphism 
$D^p \stackrel{{\cal{D}}}{\longrightarrow} D^m \rightarrow 0$. The case of a formally surjective operator, like the $div$ operator, described by the exact sequence $E \stackrel{{\cal{D}}}{\longrightarrow} F \rightarrow 0$ is now providing the exact sequence of differential modules $ 0 \rightarrow D^p \stackrel{{\cal{D}}}{\longrightarrow} D^m \rightarrow M  \rightarrow 0$ because ${\cal{D}}$ has no CC.  \\

\noindent
{\bf THEOREM 3.11}: ({\it Double Duality Test}) The procedure has $5$ steps in the operator language:  \\
$\bullet$  STEP $1$: Start with the given operator ${\cal{D}}_1$ and the corresponding differential module $M_1$.  \\
$\bullet$  STEP $2$: Construct the operator $ad({\cal{D}}_1)$.   \\
$\bullet$  STEP $3$: As any operator is the adjoint of an operator, denote by $ad({\cal{D}})$ its generating CC.  \\
$\bullet$  STEP $4$: Construct ${\cal{D}}=ad(ad({\cal{D}}))$.            \\
$\bullet$   STEP $5$: Construct the generating CC ${\cal{D}}'_1$ of ${\cal{D}}$ and compare to 
${\cal{D}}_1$.  \\
If ${\cal{D}}_1$ generates the CC of ${\cal{D}}$, we have obtained a parametrization. Otherwise, $M_1$ is not torsion-free and any new CC provides an element of $t(M_1)$.  \\

If $N_1$ is the differential module defined by $ad({\cal{D}}_1)$, it follows from the last step that $t(M_1)={ext}^1_D(N_1,D)={ext}^1(N_1)$. More generally, we have 
(See [28] p 218 for details): \\

\noindent
{\bf COROLLARY 3.12}: If $M$ is the differential module defined by {\it any} operator ${\cal{D}}$ and $N$ is the corresponding differential module defined by $ad({\cal{D}})$, then we have $t(M) = {ext}^1(N)$ with a slight abuse of language because of Definition 3.4.   \\ 

\noindent
{\bf DEFINITION 3.13}: A parametrization is said to be "{\it minimum}" if the differential module defined by ${\cal{D}}$ has a vanishing differential rank and is thus a torsion 
module ([42, 50]).  \\

\noindent
{\bf COROLLARY  3.14}: ({\it Minimum parametrization}) The procedure has $4$ steps in the operator language:  \\
$\bullet$  STEP $1$: Start with the formally exact {\it parametrizing sequence} already constructed by differential double duality.  We have thus $im({\cal{D}})=ker({\cal{D}}_1)$ and the corresponding differential module $M_1$ defined by ${\cal{D}}_1$ is torsion-free by assumption.  \\
$\bullet$  STEP $2$: Construct the adjoint sequence which is also formally exact by assumption.  \\
$\bullet$  STEP $3$: Find a {\it maximum} set of differentially independent CC $ad({\cal{D}}'):\mu \rightarrow {\nu}' $ among the generating CC $ad({\cal{D}} ):\mu \rightarrow  \nu$ of $ad({\cal{D}}_1)$ in such a way that $im(ad({\cal{D}}'))$ is a maximum free differential submodule of $im(ad({\cal{D}}))$ that is any element in $im(ad({\cal{D}}))$ is differentially algebraic over $im(ad({\cal{D}}'))$. \\
$\bullet$  STEP $4$:  Using differential duality, construct ${\cal{D}}'=ad(ad({\cal{D}}'))$.   \\
There may be many different minimum parametrizations (See [47] and section 4 for examples). \\

\noindent
{\bf EXAMPLE 3.15}: If ${\cal{D}}: \xi \rightarrow (d_{22} \xi = {\eta}^2, d_{12} \xi = {\eta}^1) $ we have ${\cal{D}}_1= ({\eta}^1, {\eta}^2) \rightarrow d_1 {\eta}^2 - d_2 {\eta}^1 = \zeta $ and the only first order generating CC of $ad({\cal{D}}_1) : \lambda \rightarrow  ( d_2 \lambda = {\mu}^1,  - d_1 \lambda = {\mu}^2) $ is $ d_1 {\mu}^1 + d_2 {\mu}^2= {\nu}' $ while $ad({\cal{D}}): ({\mu}^1, {\mu}^2) \rightarrow d_{12}{\mu}^1 + d_{22} {\mu}^2= \nu = d_2 {\nu}'$ is a second order operator like ${\cal{D}}$. \\

\noindent
{\bf EXAMPLE 3.16}: Many other examples can be found in ordinary differential control theory because it is known that a linear control system is controllable if and only if it is parametrizable (See [27, 28] for more details and examples). In our opinion, the best and simplest one is provided by the so-called double pendulum in which a rigid bar is able to move horizontally with reference position $x$ and we attach two pendula with respective length $l_1$ and $l_2$ making the (small) angles ${\theta}_1$ and ${\theta}_2$ with the vertical, the corresponding control system does not depend on the mass of each pendulum and the operatr ${\cal{D}}_1$ is defined as follows from Newton law with $d=d/dt$: 
\[     \fbox{  $  d^2x +l_1d^2{\theta}^1 +g {\theta}^1=0, \hspace{1cm}  d^2x + l_2 d^2{\theta}^2 + g{\theta}^2=0 $  }    \]
where $g$ is the gravity. The standard way used by any student of the control community, is to prove that this control system is controllable if and only if $l_1 \neq  l_2$ through a tedious computation based on  the standard Kalman test. We let the reader prove this result as an exercise and apply the previous theorem in order to work out the parametrizing operator ${\cal{D}}$ of order 4, namely:   \\
\[ \fbox{  $ \begin{array}{rcl}
  - l_1 l_2 d^4\phi - g(l_1+l_2)d^2 \phi - g^2 \phi & = & x  \\
  l_2 d^4\phi +g d^2 \phi & = & {\theta}_1  \\
  l_1 d^4 \phi + gd^2 \phi  & = & {\theta}_2
  \end{array}  $  } \]  
  The main problem is that ${\cal{D}}_1$ is trivially involutive but that {\it its adjoint is far from being even FI and the search for a Pommaret basis is quite delicate}. Indeed, multiplying the first OD equation by ${\lambda}^1$, the second by ${\lambda}^2$, adding and integrating by parts, we get: 
  \[  {\theta}^1 \rightarrow   l_1 d^2 {\lambda}^1 + g {\lambda}^1= {\mu}^1, \,\,\, {\theta}^2 \rightarrow l_2 d^2 {\lambda}^2 +  g {\lambda}^2 = {\mu}^2 , \,\,\,  
  x \rightarrow d^2 {\lambda}^1 + d^2 {\lambda}^2= {\mu}^3 \]
Multiplying the first equation by $l_2$, the second by $l_1$ and adding while taking into account the third equation, we get an equation of the form:  
\[  \fbox{ $    (l_2 {\lambda}^1 + l_1 {\lambda}^2) =\frac{1}{g}( l_2 {\mu}^1 + l_1 {\mu}^2 - l_1 l_2 {\mu}^3) = A(\mu) \in j_0(\mu) $ } \]
 Differentiating twice this equation while using the first and second equations, we also obtain: 
\[  \fbox{ $  (\frac{l_2}{l_1} {\lambda}^1 + \frac{l_1}{l_2} {\lambda}^2) =  B(\mu) \in j_2(\mu)  \,\, \Rightarrow \,\,  (l_1 - l_2) \lambda \in j_2 (\mu) $  }  \] 
When $l_1 - l_2 \neq 0$, it follows that $\lambda \in j_2(\mu)$ and, substituting in the third equation, we find the fourth order CC operator:
\[ \fbox{  $  (l_2 d^4 + g d^2){\mu}^1 + (l_1 d^4 + g d^2 ) {\mu}^2 - (l_1 l_2 d^4 + g(l_1 + l_2) d^2 + g^2 ) {\mu}^3 = \nu  $  }  \]
Multiplying by a test function $\phi$ and integrating by parts we obtain the desired parametrization. Accordingly, the differential module $M_1$ is torsion-free if and only if 
$l_1 \neq l_2$. It follows that the controllability of a control system is a {\it built in} property not depending on the choice of the single input (for example $x$) and the two outputs (for example ${\theta}^1$ and ${\theta}^2$) contrary to what the control community is still believing !. Of course, if $l_1 = l_2 = l$, setting $\theta = {\theta}^1 - {\theta}^2$, we obtain by 
subtraction $l d^2 \theta  + g \theta = 0$ and $\theta$ is a torsion element as can be seen by any reader doing the experiment. One must finally notice that the control system is controllable if and only if the adjoint of the system operator is injective (See [28] p 204-205 for details when $n=1$). With more details, if a control system is written in the Kalman {\it input/output} form $ - d y +  A y + B u = 0$, multiplying on the left by a test row vector $\lambda$ and integrating by parts, the adjoint system becomes 
$ d \lambda + \lambda A=0, \lambda B=0 \Rightarrow d \lambda B=0 \Rightarrow \lambda A B=0 \Rightarrow  \lambda A^2 B=0$ and so on. This result is showing that the Kalman controllability test amounts to the injectivity of the adjoint of the control operator {\it in the Kalman form}, {\it without any need to bring it back to a first order system}. \\
  
\noindent  
{\bf EXAMPLE 3.17} : A less academic but much more important example is the problem of parametrizing the Einstein equations, solved negatively in  $1995$ ([25]). The following 
diagram proves that {\it Einstein equations cannot be parametrized} and we shall give details in the next section:   \\
\[  \begin{array}{rcccccccl}
 & &  & 10 &\stackrel{Riemann}{\longrightarrow}  & 20 & \stackrel{Bianchi}{\longrightarrow}  & 10 & \\
 & & & \parallel & \nearrow & \downarrow &  & \downarrow &  \\
 & 4 &  \stackrel{Killing}{\longrightarrow} & 10 & \stackrel{Einstein}{\longrightarrow} & 10  & \stackrel{div}{\longrightarrow} &  4  & \longrightarrow 0 \\
  & & &  & &  &  & &  \\
 0 \longrightarrow & 4 & \stackrel{Cauchy}{\longleftarrow} & 10 & \stackrel{Einstein}{\longleftarrow} & 10 &   & &  \\
  & & & & & & & & 
\end{array}  \]
It is essential to notice that the {\it Cauchy} and {\it Killing} operators ({\it left side}) have {\it strictly nothing to do} with the {\it Bianchi} and thus {\it div} operators ({\it right side}). According to the last Corollary, the $20 - 10 = 10$ new CC are generating the torsion submodule of the differential module defined by the Einstein operator. In the last section we shall explain why such a basis of the torsion module is made by the $10$ independent components of the Weyl tensor, {\it each one killed by the Dalembertian}, a result leading to the so-called {\it Lichnerowicz waves} (in France !) [27, 28, 32, 47, 49]. \\

\noindent
{\bf EXAMPLE 3.18}: In continuum mechanics, the Cauchy stress tensor may not be symmetric in the so-called {\it Cosserat media} where the Cauchy stress equations are replaced by the Cosserat couple-stress equations which are nothing else than the adjoint of the first Spencer operator, {\it totally different} from the third [6, 26, 29]. When $n=2$, we shall see that the single Airy function has strictly nothing to do with any perturbation of the metric having three components. \\

\noindent
{\bf EXAMPLE 3.19}: A similar comment can be done for electromagnetism through the exterior derivative as the first set of Maxwell equations can be parametrized by the EM potential 1-form while the second set of Maxwell equations, {\it adjoint of this parametrization},  can be parametrized by the EM pseudo-potential [23, 24, 35]. These results are even strengthening the comments we shall make in section $4$ on the origin and existence of gravitational waves [47, 48, 49].  \\   

As a byproduct of the preceding examples, it is clear that an operator can be FI or involutive but that its adjoint may be neither involutive, nor even FI and the situation is more delicate when using double duality because each step may be as delicate as the previous one, a fact showing out the importance of the intrinsic definition of Pommaret bases that we have given and illustrated.  \\   \\

\noindent
{\bf 4) EINSTEIN EQUATIONS}  \\

Linearizing the {\it Ricci} tensor ${\rho}_{ij}$ over the Minkowski metric $\omega$, we obtain the usual second order homogeneous {\it Ricci} operator $\Omega \rightarrow R$ with $4$ terms ([7] is a fine reference):  \\
\[  2 R_{ij}= {\omega}^{rs}(d_{rs}{\Omega}_{ij}+ d_{ij}{\Omega}_{rs} - d_{ri}{\Omega}_{sj} - d_{sj}{\Omega}_{ri})= 2R_{ji}  \]
\[ tr(R)= {\omega}^{ij}R_{ij}={\omega}^{ij}d_{ij}tr(\Omega)-{\omega}^{ru}{\omega}^{sv}d_{rs}{\Omega}_{uv}  \]
We may define the $Einstein$ operator by setting $E_{ij}=R_{ij} - \frac{1}{2} {\omega}_{ij}tr(R)$ and obtain the $6$ terms [9]:  \\
\[ 2E_{ij}= {\omega}^{rs}(d_{rs}{\Omega}_{ij} + d_{ij}{\Omega}_{rs} - d_{ri}{\Omega}_{sj} - d_{sj}{\Omega}_{ri})
- {\omega}_{ij}({\omega}^{rs}{\omega}^{uv}d_{rs}{\Omega}_{uv} - {\omega}^{ru}{\omega}^{sv}d_{rs}{\Omega}_{uv})  \]
We have the (locally exact) differential sequence of operators acting on sections of vector bundles where the order of an operator is written under its arrow:  \\
\[    T \underset 1{\stackrel{Killing}{\longrightarrow}} S_2T^* \underset 2{\stackrel{Riemann}{\longrightarrow}} F_1 \underset 1{\stackrel{Bianchi}{\longrightarrow}} F_2  \]
\[    n \stackrel{{\cal{D}}}{ \longrightarrow} n(n+1)/2 \stackrel{{\cal{D}}_1}{\longrightarrow} n^2(n^2-1)/12 \stackrel{{\cal{D}}_2}{\longrightarrow} n^2(n^2-1)(n-2)/24  \]
Our purpose is now to study the differential sequence onto which its right part is projecting:  \\
\[      S_2T^* \underset 2 {\stackrel{Einstein}{\longrightarrow}} S_2T^* \underset 1{\stackrel{div}{\longrightarrow}} T^*  \rightarrow 0  \]
\[  n(n+1)/2 \longrightarrow  n(n+1)/2  \longrightarrow n  \rightarrow 0  \]
and the following adjoint sequence where we have set [35 - 38, 43] :    \\
\[   ad(T) \stackrel{Cauchy}{\longleftarrow} ad(S_2T^*)  \stackrel{Beltrami}{\longleftarrow} ad(F_1) \stackrel{Lanczos}{\longleftarrow} ad(F_2)     \]
In this sequence, if $E$ is a vector bundle over the ground manifold $X$ with dimension $n$, we may introduce the new vector bundle $ad(E)={\wedge}^nT^* \otimes E^*$ where $E^*$ is obtained from $E$ 
by inverting the transition rules exactly like $T^*$ is obtained from $T$. We have for example $ad(T)={\wedge}^nT^*\otimes T^*\simeq {\wedge}^nT^*\otimes T \simeq {\wedge}^{n-1}T^*$ because $T^*$ is isomorphic to $T$ by using the metric $\omega$.
The $10 \times 10 $ $ Einstein$ operator matrix is induced from the $10 \times 20$ $Riemann$ operator matrix and the $10 \times 4$ $div$ operator matrix is induced from the $20 \times 20$ $Bianchi$ operator matrix. We advise the reader not familiar with the formal theory of systems or operators to follow the computation in dimension $n=2$ with the $1 \times 3$ $Airy$ operator matrix, which is the formal adjoint of the $3 \times 1$ $Riemann$ operator matrix, and $n=3$ with the $6 \times 6$ $Beltrami$ operator matrix which is the formal adjoint of the $6 \times 6$ $Riemann$ operator matrix  which is easily seen to be self-adjoint up to a change of basis.\\ 
With more details, we have:  \\

\noindent
$\bullet \hspace{3mm}n=2$: The stress equations become $ d_1{\sigma}^{11}+ d_2{\sigma}^{12}=0, d_1{\sigma}^{21}+ d_2{\sigma}^{22}=0$. Their second order parametrization ${\sigma}^{11}= d_{22}\phi, {\sigma}^{12}={\sigma}^{21}= - d_{12}\phi, {\sigma}^{22}= d_{11}\phi$ has been provided by George Biddell Airy in 1863 [1] and is well known in plane elasticity [27, 35]. We get the second order system:  \\
\[ \left\{  \begin{array}{rll}
{\sigma}^{11} & \equiv d_{22}\phi =0 \\
-{\sigma}^{12} & \equiv d_{12}\phi =0 \\
{\sigma}^{22} & \equiv d_{11}\phi=0
\end{array}
\right. \fbox{ $ \begin{array}{ll}
1 & 2   \\
1 & \bullet \\  
1 & \bullet  
\end{array} $ } \]
which is involutive with one equation of class $2$, $2$ equations of class $1$ and it is easy to check that the $2$ corresponding first order CC are just the $Cauchy$  equations. Of course, the $Airy$ function ($1$ term) has absolutely nothing to do with the perturbation of the metric ($3$ terms). With more details, when $\omega$ is the Euclidean metric, we may consider the only component:   \\
\[  \begin{array}{rcl}
tr(R)  &  = &  (d_{11} + d_{22})({\Omega}_{11} + {\Omega}_{22}) - (d_{11}{\Omega}_{11} + 2 d_{12}{\Omega}_{12}+ d_{22}{\Omega}_{22})  \\
      &  =  &  d_{22}{\Omega}_{11} + d_{11}{\Omega}_{22} - 2 d_{12}{\Omega}_{12}
      \end{array}  \]
Multiplying by the Airy function $\phi$ and integrating by parts, we discover that:   \\
      \[  \fbox{ $  Airy=ad(Riemann) \,\,\,\,  \Leftrightarrow  \,\,\,\,  Riemann = ad(Airy)  $ } \]  \\
in the following adjoint differential sequences:
 \[ \fbox{  $ \begin{array}{rcccccccl}
    &  &  2  &  \underset 1 {\stackrel{Killing}{\longrightarrow}} & 3 &  \underset 2 {\stackrel{Riemann}{\longrightarrow}} &1 & \longrightarrow  &   0  \\
     &  &  &  &  &  &  &  &  \\
  0  & \longleftarrow   & 2 & \underset 1 {\stackrel{Cauchy}{\longleftarrow}}  &3  & \underset 2 {\stackrel{Airy}{\longleftarrow}}  & 1 & &   
      \end{array}  $  }  \]
      
\noindent
$\bullet \hspace{3mm} n=3$: It is quite more delicate to parametrize the $3$ PD equations: \\
\[ d_1{\sigma}^{11}+ d_2{\sigma}^{12}+ d_3{\sigma}^{13}=0,\hspace{3mm} d_1{\sigma}^{21}+ d_2{\sigma}^{22}+ d_3{\sigma}^{23}=0, \hspace{3mm} d_1{\sigma}^{31}+ d_2{\sigma}^{32}+ d_3{\sigma}^{33}=0 \]
A direct computational approach has been provided by Eugenio Beltrami in 1892 [4], James Clerk Maxwell in 1870 [19] and Giacinto Morera in 1892 [35] by introducing the $6$ {\it stress functions} ${\phi}_{ij}={\phi}_{ji}$ in the {\it Beltrami parametrization}. The corresponding system:\\
\[   \left\{  \begin{array}{rll}
{\sigma}^{11} \equiv & d_{33}{\phi}_{22}+ d_{22}{\phi}_{33}-2 d_{23}{\phi}_{23}=0  \\
-{\sigma}^{12}\equiv & d_{33}{\phi}_{12}+ d_{12}{\phi}_{33}- d_{13}{\phi}_{23}- d_{23}{\phi}_{13}=0  \\
 {\sigma}^{22}\equiv & d_{33}{\phi}_{11}+ d_{11}{\phi}_{33}-2 d_{13}{\phi}_{13}=0  \\
{\sigma}^{13}\equiv & d_{23}{\phi}_{12}+ d_{12}{\phi}_{23}- d_{22}{\phi}_{13}- d_{13}{\phi}_{22} =0 \\
-{\sigma}^{23}\equiv & d_{23}{\phi}_{11}+ d_{11}{\phi}_{23}- d_{12}{\phi}_{13}- d_{13}{\phi}_{12} =0 \\
{\sigma}^{33}\equiv & d_{22}{\phi}_{11}+ d_{11}{\phi}_{22}-2 d_{12}{\phi}_{12}=0
\end{array}
\right. \fbox{ $ \begin{array}{lll}
1 & 2 & 3   \\
1 & 2 & 3  \\
1 & 2 & 3  \\
1 & 2 &  \bullet  \\
1 & 2 & \bullet  \\
1 & 2 & \bullet
\end{array} $ } \]
is involutive with $3$ equations of class $3$, $3$ equations of class $2$ and no equation of class $1$. The three characters are thus ${\alpha}^3_2=1\times 6 - 3=3< {\alpha}^2_2=2\times 6 -3=9 < {\alpha}^1_2= 3\times 6 - 0 = 18$ and we have $dim(g_2) = {\alpha}^1_2 + {\alpha}^2_2 + {\alpha}^3_2= 18 + 9 + 3 = 30 = dim(S_2T^*\otimes S_2T^*) - dim(S_2T^*) = (6\times 6) - 6$. The $3$ CC are describing the stress equations which admit therefore a parametrization ... but without any geometric framework, in particular without any possibility to imagine that the above second order operator is {\it nothing else but} the {\it formal adjoint} of the {\it Riemann operator}, namely the (linearized) Riemann tensor with $n^2(n^2-1)/2=6$ independent components when $n=3$ [35]. \\
Breaking the canonical form of the six equations which is associated with the Janet tabular, we may rewrite the Beltrami parametrization of the Cauchy stress equations as follows, after exchanging the third row with the fourth row, keeping the ordering $\{(11)<(12)<(13)<(22)<(23)<(33)\}$:  \\
\[      \left(  \begin{array}{cccccc}
d_1& d_2 & d_3 &0 & 0 & 0 \\
 0 & d_1 &  0 & d_2 & d_3 & 0 \\
 0 & 0 & d_1 & 0 & d_2 & d_3 
\end{array}  \right)  
 \left(  \begin{array}{cccccc}
 0 & 0 & 0 & d_{33} & - 2d_{23} & d_{22} \\
 0 & - d_{33} & d_{23} & 0 & d_{13} & - d_{12}  \\
 0 & d_{23} & - d_{22} & - d_{13} & d_{12} & 0 \\
 d_{33}& 0 & - 2 d_{13} & 0 & 0 & d_{11}  \\
 - d_{23} & d_{13} & d_{12}& 0 & - d_{11} & 0 \\
 d_{22} & - 2 d_{12} & 0 & d_{11}& 0 & 0 
 \end{array} \right)  \equiv   0    \]
 as an identity where $0$ on the right denotes the zero operator. However, if  $\Omega$ is a perturbation of the metric $\omega$, the standard implicit summation used in continuum mechanics is, when $n=3$:  \\
 \[   \begin{array}{rcl}
{\sigma}^{ij}{\Omega}_{ij} & = & {\sigma}^{11}{\Omega}_{11} + 2 {\sigma}^{12}{\Omega}_{12} + 2 {\sigma}^{13}{\Omega}_{13} + {\sigma}^{22} {\Omega}_{22} + 2{\sigma}^{23}{\Omega}_{23} + {\sigma}^{33}{\Omega}_{33}  \\
   &  =  & {\Omega}_{22}d_{33}{\phi}_{11}+ {\Omega}_{33}d_{22}{\phi}_{11}- 2 {\Omega}_{23}d_{23}{\phi}_{11}+ ... \\
   &    & + {\Omega}_{23}d_{13}{\phi}_{12}+{\Omega}_{13}d_{23}{\phi}_{12}- {\Omega}_{12}d_{33}{\phi}_{12}- {\Omega}_{33}d_{12}{\phi}_{12} + ...
\end{array}  \]
because {\it the stress tensor density $\sigma$ is supposed to be symmetric}. Integrating by parts in order to construct the adjoint operator, we get:  \\
\[ \begin{array}{rcl}
 {\phi}_{11} &  \longrightarrow  &  d_{33}{\Omega}_{22} + d_{22}{\Omega}_{33} - 2 d_{23}{\Omega}_{23} \\
 {\phi}_{12} &  \longrightarrow   &  d_{13}{\Omega}_{23}+d_{23}{\Omega}_{13}-d_{33}{\Omega}_{12} - d_{12}{\Omega}_{33}
 \end{array}  \]
and so on. The identifications $ Beltrami = ad(Riemann), \,\, Lanczos = ad(Bianchi)$ in the diagram:  \\
   \[ \fbox{  $  \begin{array}{rcccccccl}  
   & 3  & \underset 1 {\stackrel{Killing}{\longrightarrow}} & 6  & \underset 2 {\stackrel{Riemann}{\longrightarrow}}& 6  & \underset 1 {\stackrel{Bianchi}{\longrightarrow }}  & 3 & \longrightarrow 0  \\
  0 \longleftarrow  & 3  & \underset 1 {\stackrel{Cauchy}{\longleftarrow}} & 6   & \underset 2 {\stackrel{Beltrami}{\longleftarrow}} 
  & 6 &  \underset 1 {\stackrel{Lanczos}{\longleftarrow}} & 3 &    
\end{array}. $   } \]
prove that the $Cauchy$ operator has nothing to do with the $Bianchi$ operator [35, 42, 47]. \\
When $\omega$ is the Euclidean metric, the link between the two sequences is established by means of the elastic constitutive relations $2 {\sigma}_{ij} = \lambda tr(\Omega) {\omega}_{ij} + 2 \mu {\Omega}_{ij}$ with the Lam\'{e} elastic constants $(\lambda, \mu) $ but mechanicians are usually setting $ {\Omega}_{ij} = 2 {\epsilon}_{ij} $. Using the standard Helmholtz decomposition ${\vec{\xi}} = {\vec{\nabla} }\varphi + {\vec{\nabla}} \wedge {\vec{\psi}} $ and substituting in the dynamical equation $ d_i {\sigma}^{ij} =  \rho d^2/dt^2 {\xi}^j $ where $\rho$ is the mass per unit volume, we get the longitudinal and transverse wave equations, namely $ \Delta \varphi - \frac{\rho}{\lambda + 2 \mu} \frac{d^2}{dt^2} \varphi = 0 $ and 
$ \Delta {\vec{\psi}} - \frac{\rho}{\mu} \frac{d^2}{dt^2} {\vec{\psi}} = 0 $, responsible for earthquakes !. \\

Then, taking into account the factor $2$ involved by multiplying the second, third and fifth row by $2$, we get the new $6\times 6$ operator matrix with rank $3$ which is clearly self-adjoint:  \\
\[   \left(  \begin{array}{cccccc}
 0 & 0 & 0 & d_{33} & - 2d_{23} & d_{22} \\
 0 & - 2d_{33} & 2d_{23} & 0 & 2d_{13} & - 2d_{12}  \\
 0 & 2d_{23} & - 2d_{22} & - 2d_{13} & 2d_{12} & 0 \\
 d_{33}& 0 & - 2 d_{13} & 0 & 0 & d_{11}  \\
 - 2d_{23} & 2d_{13} & 2d_{12}& 0 & - 2d_{11} & 0 \\
 d_{22} & - 2 d_{12} & 0 & d_{11}& 0 & 0 
 \end{array} \right)     \]

{\it Surprisingly}, the Maxwell parametrization is obtained by keeping ${\phi}_{11}=A, {\phi}_{22}=B, {\phi}_{33}=C$ while setting ${\phi}_{12}={\phi}_{23}={\phi}_{31}=0$ in order to obtain the system:\\
\[   \left\{  \begin{array}{rl}
{\sigma}^{11} \equiv& d_{33}B + d_{22}C=0  \\
{\sigma}^{22}\equiv & d_{33}A+ d_{11}C =0 \\
- {\sigma}^{23}\equiv & d_{23}A=0  \\
{\sigma}^{33}\equiv & d_{22}A+ d_{11}B=0  \\
- {\sigma}^{13}\equiv & d_{13}B=0 \\
- {\sigma}^{12}\equiv & d_{12}C=0
\end{array}
\right. \fbox{ $ \begin{array}{lll}
1 & 2 & 3   \\
1 & 2 & 3  \\
1 & 2 & \bullet  \\
1 & 2 &  \bullet  \\
1 & \bullet & \bullet  \\
1 & \bullet & \bullet
\end{array} $ } \]
{\it This system may not be involutive} and no CC can be found "{\it a priori} " because the coordinate system is surely not $\delta$-regular. Effecting the linear change of coordinates 
${\bar{x}}^1 = x^1, {\bar{x}}^2 = x^2, {\bar{x}}^3 = x^3 + x^2 + x^1 $ and taking out the bar for simplicity, we obtain the homogeneous involutive system with a quite tricky Pommaret basis:  \\
\[   \left\{  \begin{array}{l}
d_{33}C+ d_{13}C+ d_{23}C+ d_{12}C=0  \\
d_{33}B+ d_{13}B=0  \\
d_{33}A+ d_{23}A=0  \\
d_{23}C +d_{22}C - d_{13}C - d_{13}B - d_{12}C =0  \\
d_{23}A - d_{22}C + d_{13}B + 2 d_{12}C - d_{11}C=0  \\
d_{22}A + d_{22}C - 2 d_{12}C + d_{11}C + d_{11}B=0
\end{array} \right. 
\fbox{ $ \begin{array}{lll}
1 & 2 & 3   \\
1 & 2 & 3  \\
1 & 2 &  3  \\
1 & 2 &  \bullet  \\
1 &  2 & \bullet  \\
1 &  2 & \bullet
\end{array} $ } \]
It is easy to check that the $3$ CC obtained just amount to the desired $3$ stress equations when coming back to the original system of coordinates. However, the three characters are different as we have now ${\alpha}^3_2=3 - 3=0 < {\alpha}^2_2= 2\times 3 -3=3 < {\alpha}^1_2=3 \times 3 - 0=9$ with sum equal to $dim(g_2)=6\times 3 - 6=18 -6= 12$. {\it We have thus a minimum parametrization}. \\

Again, {\it if there is a geometrical background, this change of local coordinates is hiding it totally}. Moreover, we notice that the stress functions kept in the procedure are just the ones on which ${\partial}_{33}$ is acting. The reason for such an apparently technical choice is related to very general deep arguments in the theory of differential modules, namely the fact that one can always find a minimum parametrization, whenever {\it a } parametrization is known to exist. \\
The following new minimum parametrization does not seem to be known:  \\
\[   \left\{  \begin{array}{rll}
{\sigma}^{11} \equiv & d_{33}{\phi}_{22} = 0  \\
-{\sigma}^{12}\equiv & d_{33}{\phi}_{12} = 0  \\
 {\sigma}^{22}\equiv & d_{33}{\phi}_{11} = 0  \\
{\sigma}^{13}\equiv & d_{23}{\phi}_{12} - d_{13}{\phi}_{22} =0 \\
-{\sigma}^{23}\equiv & d_{23}{\phi}_{11} - d_{13}{\phi}_{12} =0 \\
{\sigma}^{33}\equiv & d_{22}{\phi}_{11}+ d_{11}{\phi}_{22}-2 d_{12}{\phi}_{12} = 0
\end{array}
\right. \fbox{ $ \begin{array}{lll}
1 & 2 & 3   \\
1 & 2 & 3  \\
1 & 2 & 3  \\
1 & 2 &  \bullet  \\
1 & 2 & \bullet  \\
1 & 2 & \bullet
\end{array} $ } \]

When $n=4$, taking the adjoint of the second order PD equations defining the so-called gravitational waves, we have proved in many books  ([47]) or papers ([48, 49]) the following crucial theorem which is showing that the Einstein operator is useless contrary to the classical GR literature ([32]).  \\

\noindent
{\bf THEOREM  4.1}: {\it The GW equations are defined by the adjoint of the Ricci operator which is not self-adjoint contrary to the Einstein operator which is indeed self-adjoint}.  \\

We finally prove that this result only depends on the second order jets of the conformal group of transformations of space-time, {\it a result highly not evident at first sight for sure}. \\ \\

\noindent
{\bf 5) CONFORMAL GROUP}  \\

We start proving that the structure of the conformal with $(n+1)(n+2)/2$ parameters may not be related to a classification of Lie algebras [43]). \\
For $n=1$, the simplest such group of transformations of the real line with $3$ parameters is the projective group defined by the Schwarzian third order OD equation:
\[   \Phi(y, y_x, y_{xx}, y_{xxx}) \equiv \frac{y_{xxx}}{y_x} - \frac{3}{2} (\frac{y_{xx}}{y_x})^2 = \nu(x)\]
with linearization the only third order Medolaghi equation with symbol $g_3=0$:
\[    L({\xi}_3)\nu \equiv  {\xi}_{xxx} + 2 \nu (x) {\xi}_x + \xi {\partial}_x \nu (x)= 0  \]
When $\nu = 0$, the general solution is simply $\xi = a x^2 + b x + c$ with $3$ parameters. There is no CC. \\

For $n=2$, eliminating the conformal factor in the case of the Euclidean metric of the plane provides the two Cauchy-Riemann equations defining the infinitesimal complex transformations of the plane. The {\it only possibility} coherent with homogeneity is thus to consider the following system and to prove that it is defining a  system of 
infinitesimal Lie equations, leading to $6$  infinitesimal generators, namely: {\it 2 translations + 1 rotation + 1 dilatation + 2 elations}:  \\
\[  \left \{ \begin{array}{c}
  {\xi}^k_{ijr}=0   \\
   {\xi}^2_{22} - {\xi}^1_{12}=0, {\xi}^1_{22} + {\xi}^2_{12}=0, {\xi}^2_{12} - {\xi}^1_{11}=0, 
{\xi}^1_{12} + {\xi}^2_{11}=0   \\
  {\xi}^2_2 - {\xi}^1_1=0, {\xi}^1_2 + {\xi}^2_1=0   
  \end{array} \right.   \]
\[\{ {\theta}_1={\partial}_1, {\theta}_2={\partial}_2, {\theta}_3=x^1 {\partial}_2 - x^2 {\partial}_1, {\theta}_4=x^1 {\partial}_1 + x^2 {\partial}_2, {\theta}_5= - \frac{1}{2} ((x^1)^2 + (x^2)^2){\partial}_1  + x^1 (x^1 {\partial}_1 + x^2 {\partial}_2) , {\theta}_6 \}  \]
with the elation ${\theta}_6$ obtained from ${\theta}_5$ by exchanging $x^1$ with $x^2$. We have ${\hat{g}}_3=0$ when $n=1,2$.  \\

\noindent
{\bf LEMMA 5.1}: We have ([33]):  \\
$\bullet$ ${\hat{g}}_1$ is finite type with ${\hat{g}}_3=0, \forall n \geq 3$.  \\
$\bullet$ ${\hat{g}}_2$ is $2$-acyclic when $n\geq 4$.  \\
$\bullet$ ${\hat{g}}_2$ is $3$-acyclic when $n\geq 5$.  \\

In order to convince the reader that both classical and conformal differential geometry must be revisited, let us prove that the analogue of the Weyl tensor is made by a third order operator when $n=3$, a result which is neither known nor acknowledged today. We shall proceed by diagram chasing as the local computation done by using computer algebra does not provide any geometric insight (See arXiv:1603.05030 and [37] for the details). We have $E=T$ and $dim({\hat{F}}_0)=5$ in the following commutative diagram providing 
${\hat{F}}_1$:  \\

\[   \begin{array}{rcccccccccl}
  & & 0  &  & 0  & &  0  &           \\
  & & \downarrow  &  &  \downarrow & & \downarrow &    \\
0  & \rightarrow & {\hat{g}}_4 & \rightarrow &S_4T^*\otimes T& \rightarrow &S_3 T^ *\otimes {\hat{F}}_0 & \rightarrow & {\hat{F}}_1 & \rightarrow  &  0 \\
& & \downarrow  &  &  \downarrow & & \parallel & &    \\
0 & \rightarrow & T^*\otimes {\hat{g}}_3 & \rightarrow &T^*\otimes S_3T^*\otimes T& \rightarrow & T^* \otimes S_2 T^ *\otimes {\hat{F}}_0 &  \rightarrow &0 & &   \\
& & \downarrow  &  &  \downarrow & & \downarrow & & &    \\
0 & \rightarrow &{\wedge}^2 T^*\otimes {\hat{g}}_2 & \rightarrow &{\wedge}^2T^* \otimes S_2T^* \otimes T& \rightarrow & {\wedge}^2 T^* \otimes T^* \otimes {\hat{F}}_0 &\rightarrow &    0   & &    \\
& & \downarrow  &  &  \downarrow & & \downarrow & & &   \\
0 & \rightarrow &{\wedge}^3 T^*\otimes {\hat{g}}_1 & \rightarrow &{\wedge}^3T^*\otimes T^*\otimes T& \rightarrow & {\wedge}^3 T^* \otimes {\hat{F}}_0 & \rightarrow  & 0  & \\
&   & \downarrow  &  &  \downarrow & & \downarrow & & &    \\
 &&  0  &  & 0  & &  0  &  &   &        
\end{array}  \]

\[   \begin{array}{rcccccccccl}
  & &   &  & 0  & &  0  &           \\
  & &   &  &\downarrow   & & \downarrow &    \\
  &  & 0 & \rightarrow &  45  & \rightarrow &  50 & \rightarrow & 5 & \rightarrow  &  0 \\
  & &   &  &  \downarrow & & \downarrow & &    \\
  &  & 0 & \rightarrow &  90 & \rightarrow &  90 &  \rightarrow &0 & &   \\
& & \downarrow  &  &  \downarrow & & \downarrow & & &    \\
0 & \rightarrow & 9 & \rightarrow & 54 & \rightarrow & 45 &\rightarrow & 0 & &    \\
& & \downarrow  &  &  \downarrow & & \downarrow &  &  \\
0 & \rightarrow & 4  & \rightarrow & 9 & \rightarrow &  5   &\rightarrow &0& \\
&   & \downarrow  &  &  \downarrow &  &  \downarrow & & &     \\
 &&  0  &  & 0  & & 0   &  &   &        
\end{array}  \]
A delicate double circular chase provides ${\hat{F}}_1=H^2_2({\hat{g}}_1)$ in the short exact sequence: \\
\[   0 \longrightarrow {\hat{F}}_1 \longrightarrow {\wedge}^2 T^* \otimes {\hat{g}}_2 \stackrel{\delta}{\longrightarrow} {\wedge}^3 T^* \otimes {\hat{g}}_1 \rightarrow 0  \hspace{2cm}
      0 \longrightarrow 5 \longrightarrow  9  \stackrel{\delta}{\longrightarrow}  4  \rightarrow 0  \]
but we have to prove that the map $\delta$ on the right is surjective, a result that it is almost impossible to prove in local coordinates. Let us prove it by means of circular diagram chasing in the preceding commutative diagram as follows. Lift any $a \in {\wedge}^3 T^* \otimes {\hat{g}}_1 \subset {\wedge}^3 T^* \otimes T^* \otimes T$ to $b \in {\wedge}^2 T^* \otimes  S_2T^* \otimes T$ because the vertical $\delta$-sequence for $S_4 T^*$ is exact. Project it by the symbol map ${\sigma}_1({\hat{\Phi}})$ to $c \in {\wedge}^2 T^* \otimes  T^* \otimes {\hat{F}}_0$. Then, lift $c$ to $d \in T^* \otimes S_2 T^* \otimes {\hat{F}}_0$ that we may lift {\it backwards horizontally} to $e \in T^* \otimes S_2 T^* \otimes T$ to which we may apply $\delta$ to obtain $f \in {\wedge}^2 T^* \otimes S_2T^* \otimes T$. By commutativity, both $f $ {\it and} $b$ map to $c$ and the difference $f - b$ maps thus to zero. Finally, we may find $g \in {\wedge}^2 T^* \otimes {\hat{g}}_2 $ such that $ b = g + \delta (e)$ and we obtain thus $a = \delta (g) + {\delta}^2 (e) = \delta (g)$, proving therefore the desired surjectivity. \\
We have $10$ parameters: {\it 3  translations +  3 rotations + 1 dilatation + 3 elations} and the totally unexpected formally exact sequences on the jet level are thus, showing in particular that second order CC do not exist:  \\
\[  0 \rightarrow {\hat{R}}_3 \rightarrow J_3(T)  \rightarrow J_2({\hat{F}}_0) \rightarrow 0 \,\, \Rightarrow \,\, 
0 \rightarrow 10 \rightarrow 60 \rightarrow 50 \rightarrow 0  \] 
\[     0  \rightarrow {\hat{R}}_4 \rightarrow J_4(T)  \rightarrow J_3({\hat{F}}_0) \rightarrow {\hat{F}}_1 \rightarrow 0 \,\,\, \Rightarrow \,\, 
0 \rightarrow 10 \rightarrow 105 \rightarrow 100 \rightarrow 5 \rightarrow 0 \]
\[     0  \rightarrow {\hat{R}}_5 \rightarrow J_5(T)  \rightarrow J_4({\hat{F}}_0) \rightarrow J_1({\hat{F}}_1) \rightarrow {\hat{F}}_2  \rightarrow 0 \,\,\, \Rightarrow \,\, 
0 \rightarrow 10 \rightarrow 168 \rightarrow 175 \rightarrow 20 \rightarrow   3    \rightarrow 0 \]
We obtain the minimum differential sequence, {\it which is nervertheless  not a Janet sequence}:
\[   0 \rightarrow {\hat{\Theta}} \rightarrow T \underset 1{\stackrel{\hat{\cal{D}}}{\rightarrow}} {\hat{F}}_0 \underset 3{\rightarrow} {\hat{F}}_1 \underset 1{\rightarrow}  {\hat{F}}_2  \rightarrow 0 \,\, \Rightarrow \,\, 
0 \rightarrow {\hat{\Theta}} \rightarrow 3 \underset 1{\stackrel{\hat{\cal{D}}}{\rightarrow}} 5 \underset 3{\rightarrow} 5 \underset 1{\rightarrow} 3 \rightarrow 0   \]
with $\hat{\cal{D}}$ the conformal Killing operator and vanishing Euler-Poincar\'{e} characteristic  
$3 - 5 + 5 - 3 = 0$.   \\

When $n=4$, we have $15$ parameters: {\it 4 translations + 6 rotations + 1 dilatation + 4 elations}.  \\

 When $n=4$ and ${\hat{g}}_3=0  \Rightarrow {\hat{g}}_4 = 0 \Rightarrow  {\hat{g}}_5 = 0 $ in the conformal case, we have the commutative diagram with exact vertical long $\delta$-sequences {\it but the left one} and where the second row proves that there cannot exist first order Bianchi-like identities for the Weyl tensor, contrary to what is still believed today:

\[  \begin{array}{rcccccccccl}
  & &  0 & & 0 & & 0 &  & 0 &  & \\
  & & \downarrow & & \downarrow & & \downarrow & & \downarrow & & \\
0 & \rightarrow & {\hat{g}}_4 & \rightarrow &  S_4T^*\otimes T & \rightarrow & S_3T^*\otimes {\hat{F}}_0& \rightarrow & T^* \otimes {\hat{F}}_1 & \rightarrow & 0  \\
  & & \hspace{2mm}\downarrow  \delta  & & \hspace{2mm}\downarrow \delta & &\hspace{2mm} \downarrow \delta & & \parallel & & \\
0 & \rightarrow& T^*\otimes {\hat{g}}_3 &\rightarrow &T^*\otimes S_3T^*\otimes T & \rightarrow &T^*\otimes  S_2 T^*\otimes {\hat{F}}_0 &\rightarrow & T^* \otimes {\hat{F}}_1  & \rightarrow & 0 \\
  & &\hspace{2mm} \downarrow \delta &  &\hspace{2mm} \downarrow \delta & &\hspace{2mm}\downarrow \delta &  & \downarrow  & &  \\
0 & \rightarrow & {\wedge}^2T^*\otimes {\hat{g}}_2 & \rightarrow & {\wedge}^2T^*\otimes  S_2 T^*\otimes T & \rightarrow & {\wedge}^2T^*\otimes  T^* \otimes {\hat{F}}_0 & \rightarrow & 0 &&  \\
 &  &\hspace{2mm}\downarrow \delta  &  & \hspace{2mm} \downarrow \delta  &  & \hspace{2mm} \downarrow \delta & &  & & \\
0 & \rightarrow & {\wedge}^3T^*\otimes {\hat{g}}_1 & =  & {\wedge}^3T^*\otimes T^* \otimes  T  & \rightarrow   &  {\wedge}^3 T^* \otimes {\hat{F}}_0  & \rightarrow  & 0 &  & \\
 &   & \hspace{2mm} \downarrow \delta &  & \hspace{2mm}  \downarrow \delta  &  &  \downarrow    &  &  &  &\\
0  & \rightarrow  &  {\wedge}^4 T^* \otimes T  & =  & {\wedge}^4 T^* \otimes T & \rightarrow & 0 &  &  &&  \\
  &   & \downarrow &  &  \downarrow  &  &  &  &  &  &  \\
   & &  0  &  &  0  &  &  &  & &  &
\end{array}  \]

\[  \begin{array}{rcccccccccl}
  & &   & & 0 & & 0 &  & 0 & &  \\
  & & & & \downarrow  & & \downarrow & & \downarrow &  & \\
   & & 0 & \rightarrow &  140 & \rightarrow & 180 & \rightarrow & 40 & \rightarrow &  0  \\
  &  & & & \hspace{2mm} \downarrow \delta  & &\hspace{2mm} \downarrow \delta & & \parallel & &  \\
 & & 0 & \rightarrow & 320& \rightarrow & 360 &\rightarrow & 40  & \rightarrow & 0  \\
 &  & \downarrow  &  &\hspace{2mm} \downarrow \delta & &\hspace{2mm}\downarrow \delta &  & \downarrow  &  &  \\
0 & \rightarrow & 24 & \rightarrow & 240 & \rightarrow & 216 & \rightarrow & 0 & &  \\
  & &\hspace{2mm}\downarrow \delta  &  & \hspace{2mm} \downarrow \delta  &  & \hspace{2mm} \downarrow \delta  & & & & \\
0 & \rightarrow & 28 & \rightarrow  & 64 &\rightarrow   & 36    & \rightarrow  & 0 &  &  \\
  &  &  \hspace{2mm} \downarrow \delta  &  & \hspace{2mm} \downarrow \delta &  & \downarrow &  &  &  & \\
 0 & \rightarrow  &  4  & =  & 4  &\rightarrow  & 0 &  &  & &  \\
  &  &  \downarrow &  & \downarrow & &  &  &  &  &  \\
  &  &   0   &  &  0  &  &  &  &  &  &  
\end{array}  \]
A diagonal snake chase proves that ${\hat{F}}_1 \simeq H^2 ({\hat{g}}_1)$. However, we have the $\delta$-sequence:  \\
\[   0 \rightarrow  T^*\otimes  {\hat{g}}_2 \stackrel{\delta}{\rightarrow} {\wedge}^2 T^* \otimes {\hat{g}}_1 \stackrel{\delta}{\rightarrow } {\wedge}^3 T^*\otimes T \rightarrow 0  \]
We obtain $dim(B^2_2({\hat{g}}_1))= 4 \times 4 = 16$ and let the reader prove as before that the map $\delta$ on the right is surjective, a result leading to
$dim(Z^2_2({\hat{g}}_1))= (6 \times (6 + 1)) - (4 \times 4) = 42 - 16 = 26$. The Weyl tensor has thus $dim({\hat{F}}_1) = 26 - 16 = 10$ components, a way that must be compared to the standard one that can be found in the GR literature. We obtain the minimum differential sequence, {\it which is nervertheless  not a Janet sequence}:
\[   0 \rightarrow {\hat{\Theta}} \rightarrow T \underset 1{\stackrel{\hat{\cal{D}}}{\rightarrow}} {\hat{F}}_0 \underset 2{\rightarrow} {\hat{F}}_1 \underset 2{\rightarrow}  {\hat{F}}_2  \underset 1{\rightarrow}  {\hat{F}}_3  \rightarrow 0 \,\, \Rightarrow \,\, 
0 \rightarrow {\hat{\Theta}} \rightarrow 4  \underset 1{\stackrel{\hat{\cal{D}}}{\rightarrow}} 9 \underset 2{\rightarrow} 10 \underset 2{\rightarrow} 9  \underset 1{\rightarrow}  4  \rightarrow 0   \]

As a byproduct, we end this paper with the following {\it fundamental diagram} $II$ first presented in $1983$ (See [23], p 430 and the reference 87 p 560) but still not yet acknowledged 
as it only depends on the Spencer $\delta$-map, explaining both the splitting vertical sequence on the right  and the link existing between the Ricci vector bundle and the symbol bundle ${\hat{g}}_2 \simeq T^*$ of second order jets of conformal elations. Needless to say that the diagonal chase providing the isomorphism $Ricci \simeq S_2T^*$ could not be even imagined by using classical methods because its involves Spencer $\delta$-cohomology with the standard notations $B=im(\delta), Z= ker(\delta), H=Z/B$ for {\it coboundary, cocycle, cohomology} at ${\wedge}^sT^* \otimes g_{q+r}$ when $g_{q+r}$ is the $r$-prolongation of a symbol $g_q$. It is important to notice that all the bundles appearing in this diagram only depend on the metric $\omega$ but {\it not} on any conformal factor. \\
 
\[  \fbox{  $   \begin{array}{rcccccccccl}
&  &  &  &  &  &  &  &  &  & \\
  &  &  &  &  &  &   &  &  0  & & \\
  &  &  &  &  &  &   &  & \downarrow &  & \\
  &  &  &  &  &  &  0  &  &  Ricci & &  \\
  &  &  &  &  &  &  \downarrow &  & \downarrow &  &  \\
  &  &  &  &  0 & \longrightarrow & Z^2_1(g_1) & \longrightarrow & Riemann & \longrightarrow & 0  \\
  &  &  &  &   \downarrow &  & \downarrow &  & \downarrow &  &  \\
  &  & 0 & \longrightarrow & T^* \otimes {\hat{g}}_2 & \stackrel{\delta}{\longrightarrow} & Z^2_1({\hat{g}}_1) & \longrightarrow & Weyl  & \longrightarrow & 0  \\
  &  &  &  &  \downarrow &  &  \downarrow &  &  \downarrow &  &  \\
 0 &  \longrightarrow & S_2T^* &  \stackrel{\delta}{\longrightarrow} &  T^* \otimes T^* & \stackrel{\delta}{\longrightarrow} & {\wedge}^2 T^* &  \longrightarrow & 0 &  &  \\      
  &  &  &  & \downarrow &  &  \downarrow &  &  &  &  \\
  &  &  &  &  0 &  &  0  &  &  &  &\\
  &  &  &  &  &  &  &  &  &  & 
  \end{array}  $  }  \]   \\
 
\[    \begin{array}{rcccccccccl}
  &  &  &  &  &  &   &  &  0  & & \\
  &  &  &  &  &  &   &  & \downarrow &  & \\
  &  &  &  &  &  &  0  &  &  10  & &  \\
  &  &  &  &  &  &  \downarrow &  & \downarrow &  &  \\
  &  &  &  &  0 & \longrightarrow & 20 & \longrightarrow & 20 & \longrightarrow & 0  \\
  &  &  &  &   \downarrow &  & \downarrow &  & \downarrow &  &  \\
  &  & 0 & \longrightarrow &  16  & \stackrel{\delta}{\longrightarrow} & 26 & \longrightarrow & 10 & \longrightarrow & 0  \\
  &  &  &  &  \downarrow &  &  \downarrow &  &  \downarrow &  &  \\
 0 &  \longrightarrow & 10 &  \stackrel{\delta}{\longrightarrow} &  16  & \stackrel{\delta}{\longrightarrow} & 6 &  \longrightarrow & 0 &  &  \\      
  &  &  &  & \downarrow &  &  \downarrow &  &  &  &  \\
  &  &  &  &  0 &  &  0  &  &  &  &
  \end{array}  \] 
  
When $n=4$, we have explained (See [24, 41] and the recent [44]) that the splitting horizontal lower sequence provides an isomorphism $T^* \otimes {\hat{g}}_2 \simeq T^* \otimes T^* \simeq S_2T^* \oplus {\wedge}^2 T^*$ locally described by $ (R_{ij}, F_{ij})$ in which $(R_{ij})$ is the GR part and $(F_{ij})$ the EM part as a unification of gravitation and electromagnetism, only depending thus on the second order jets of conformal transformations, {\it contrary to the philosophy of GR today but in a coherent way with the dream of H. Weyl} ([55]). Introducing the Weyl algebroid ${\tilde{R}}_2$ with $11$ parameters, we finally notice that $T^* \otimes {\hat{g}}_2 =T^*\otimes ({\hat{R}}_2 / {\tilde{R}}_2) = (T^* \otimes  {\hat{R}}_2) / (T^* \otimes  {\tilde{R}}_2)= {\hat{C}}_1 / {\tilde{C}}_1$, a result contradicting the mathematical foundations of classical gauge theory while allowing to understand the confusion done by E. Cartan and followers between "curvature alone" ($F_1$) and "curvature + torsion" ($C_2$) while the EM field comes from a section of $C_1$ (See [44] for more details).  \\   \\

\noindent
{\bf 6) CONCLUSION}   \\

It is not so well known that a classical OD control system defined by a surjective operator ${\cal{D}}$ is controllable if and only if the operator $ad({\cal{D}})$ is injective or, 
{\it equivalently}, if the operator ${\cal{D}}$ can be parametrized. The simplest example is the Kalman system $y_x = A y + B u$ with input $u$ and output $y$ leading to the Kalman test but this result cannot be extended to an arbitrary PD control system with many independent variables ([27]). In this case, one needs the double differential duality test for checking if the corresponding differential module $M$ is torsion-free or, {\it equivalently}, if ${\cal{D}}$ can be parametrized. Then, one has in general to use twice the PP procedure which is already delicate for the OD case (See the double pendulum) but may become awful for the PD case, like in the study of the Killing operator for the Kerr metric ([45]). However, it is a fact that both the control, computer algebra, physics communities largely refused to use the Spencer operator and we don't speak about the mechanical community still not accepting that the Cosserat couple-stress equations are nothing else than the adjoint of the first Spencer operator in the Spencer sequence for the group of rigid motions in space ([26, 29]) or that the second set of Maxwell equations in electromagnetism are similarly induced by the adjoint of the first Spencer operator for the group of conformal transformations in space-time along the dream of Weyl ([55]).  We have tried to explain in this paper what are the negative consequences on the origin and existence of gravitational waves ([47]). As a conclusion, we can only hope that such a poor effective situation will indeed be improved in the future.  \\   \\

\noindent
{\bf REFERENCES}  \\

\noindent
[1] Airy, G.B.:  On the Strains in the Interior of Beams, Phil. Trans. Roy. Soc. London, 153 (1863) 49-80.  \\  
\noindent
[2] Albert, M., Fetzer, M., Saenz-de-Cabezon, E., Seiler, W.M., On the Free Resolutions Induced by a Pommaret Basis, Journal of Symbolic Computation, 68 (2015) 4-26.  \\
https://doi.org/10.1016/j.jsc.2014.09.008  \\
\noindent
[3] Apel, J.: The Theory of Involutive Divisions and an Application to Hilbert Functions Computations, J. Symbolic Computations, 25 (1998) 683-704.  \\
\noindent
[4] Beltrami, E.: Osservazioni sulla Nota Precedente, Atti della Accademia Nazionale dei Lincei Rend., 1, 5 (1892) 141-142; Collected Works, t IV . \\
\noindent
[5] Binaei, B., Hashemi, A., Seiler, W.M.: Computations of Pommaret Bases Using Syzygies, arXiv:1809.10971  \\
\noindent
[6] Cosserat, E., Cosserat, F.: Th\'{e}orie des Corps D\'{e}formables, Hermann, Paris (1909).\\
\noindent
[7] Foster, J., Nightingale, J.D.: A Short Course in General Relativity, Longman (1979).  \\
\noindent
[8] Gerdt, V. P.: Gr\"{o}bner bases and involutive methods for algebraic and differential equations, , A. Yu.Mathematical and computer modelling, vol. 25, no. 8-9, pp. 75-90, 1997. DOI: 10.1016/S0895-7177(97)00060-5.   \\
\noindent
[9] Gerdt, V.P.: Involutive Division Technique: Some Generalizations and Optimizations, Zapiski Nauchnykh Seminarov POMI (St.Petersburg) 258 (1999) 185-206. To appear in J. Math. Sci. 258 (2000).  \\
\noindent
[10] Gerdt, V.P.: On the Relation Between Pommaret and Janet Bases, arXiv:math/0004100.  \\
\noindent
[11] Gerdt, V.P., Blinkov, Yu.A.: Involutive Bases of Polynomial Ideals. Math. Comp. Sim. 45 (1998) 519-542.   \\
\noindent
[12] Gerdt, V.P., Blinkov, Yu.A.: Minimal Involutive Bases. Math. Comp. Simul. 45 (1998) 543-560.  \\
\noindent
[13] Gerdt, V.P., Zinin, M.V.: A Pommaret Division Algorithm for Computing Gr\"{o}bner Bases in Boolean rings, ISSAC'08, July 20-23, 2008, Hagenberg, Austria.  \\
\noindent
[14] Goldschmidt, H.: Prolongations of Linear Partial Differential Equations: I Inhomogeneous equations, Ann. Scient. Ec. Norm. Sup., 4 (1968) 617-625.  https://doi.org/10.24033/asens.1173\\
\noindent
[15] Janet, M.: Sur les Syst\`{e}mes aux D\'{e}riv\'{e}es Partielles, Journal de Math., 8 (1920) 65-151. \\
\noindent 
[16] Kashiwara, M.: Algebraic Study of Systems of Partial Differential Equations, M\'{e}moires de la Soci\'{e}t\'{e} Math\'{e}matique de France, 63 (1995) (Transl. from Japanese of his 1970 Master Thesis).  \\
\noindent
[17] Kumpera, A., Spencer, D.C.: Lie Equations, Ann. Math. Studies 73, Princeton University Press, Princeton (1972).\\
\noindent
[18] Macaulay, F.S.: The Algebraic Theory of Modular Systems, Cambridge Tract 19, Cambridge University Press, London, 1916 (Reprinted by Stechert-Hafner Service Agency, New York, 1964).  \\
\noindent
[19] Maxwell, J.C.: On Reciprocal Figures, Frames and Diagrams of Forces, Trans. Roy. Soc. Ediinburgh, 26 (1870) 1-40.  \\
\noindent
[20] Northcott, D.G.: An Introduction to Homological Algebra, Cambridge university Press (1966).  \\
\noindent
[21] Pommaret, J.-F.: Systems of Partial Differential Equations and Lie Pseudogroups, Gordon and Breach, New York (1978); Russian translation: MIR, Moscow,(1983).\\
\noindent
[22] Pommaret, J.-F.: Differential Galois Theory, Gordon and Breach, New York (1983).\\
\noindent
[23] Pommaret, J.-F.: Lie Pseudogroups and Mechanics, Gordon and Breach, New York (1988).\\
\noindent
[24] Pommaret, J.-F.: Partial Differential Equations and Group Theory, Kluwer (1994).\\
http://dx.doi.org/10.1007/978-94-017-2539-2    \\
\noindent
[25] Pommaret, J.-F.: Dualit\'{e} Diff\'{e}rentielle et Applications, Comptes Rendus Acad\'{e}mie des Sciences Paris, S\'{e}rie I, 320 (1995) 1225-1230.  \\
\noindent
[26] Pommaret, J.-F.: Fran\c{c}ois Cosserat and the Secret of the Mathematical Theory of Elasticity, Annales des Ponts et Chauss\'ees, 82 (1997) 59-66 (Translation by D.H. Delphenich).  \\
\noindent
[27] Pommaret, J.-F.: Partial Differential Control Theory, Kluwer, Dordrecht (2001) (Zbl 1079.93001).   \\
\noindent
[28] Pommaret, J.-F.: Algebraic Analysis of Control Systems Defined by Partial Differential Equations, in "Advanced Topics in Control Systems Theory", Springer, Lecture Notes in Control and Information Sciences 311 (2005) Chapter 5, pp. 155-223.\\
\noindent
[29] Pommaret, J.-F.: Parametrization of Cosserat Equations, Acta Mechanica, 215 (2010) 43-55.\\
http://dx.doi.org/10.1007/s00707-010-0292-y  \\
\noindent
[30] Pommaret, J.-F.: Macaulay Inverse Systems Revisited, Journal of symbolic Computation, 46 (2011) 1049-1069. https://doi.org/10.1016/j.jsc.2011.05.007  \\
\noindent
[31] Pommaret, J.-F.: Spencer Operator and Applications: From Continuum Mechanics to Mathematical Physics, in "Continuum Mechanics-Progress in Fundamentals and Engineering Applications", Dr. Yong Gan (Ed.), ISBN: 978-953-51-0447--6, InTech (2012) Available from: \\
http://dx.doi.org/10.5772/35607   \\
\noindent
[32] Pommaret, J.-F.: The Mathematical Foundations of General Relativity Revisited, Journal of Modern Physics, 4 (2013) 223-239. https://doi.org/10.4236/jmp.2013.48A022   \\
 \noindent
[33] Pommaret, J.-F.: From Thermodynamics to Gauge Theory: The Virial Theorem Revisited. In: Gauge Theories and Differential Geometry, NOVA Science Publishers, Chapter 1 (2015) 1-44. https://arxiv.org/abs/1504.04118  \\
\noindent
[34] Pommaret, J.-F.: Relative Parametrization of Linear Multidimensional Systems, Multidim. Syst. Sign. Process., 26 (2015) 405-437. https://doi.org/10.1007/s11045-013-0265-0 \\
\noindent
[35] Pommaret, J.-F.: Airy, Beltrami, Maxwell, Einstein and Lanczos Potentials Revisited, Journal of Modern Physics, 7 (2016) 699-728. https://dx.doi.org/10.4236/jmp.2016.77068   \\
\noindent
[36] Pommaret, J.-F.: Why Gravitational Waves Cannot Exist, Journal of Modern Physics, 8 (2017) 2122-2158. https://doi.org/104236/jmp.2017.813130  (https://arxiv.org/abs/1708.06575  \\
\noindent
[37] Pommaret, J.-F.: Deformation Theory of Algebraic and Geometric Structures, Lambert Academic Publisher (LAP), Saarbrucken, Germany (2016). http://arxiv.org/abs/1207.1964 \\
\noindent
[38] Pommaret, J.-F.: New Mathematical Methods for Physics, Mathematical Physics Books, Nova Science Publishers, New York (2018) 150 pp.  \\
\noindent
[39] Pommaret, J.-F.: Differential Homological Algebra and General Relativity. Journal of Modern Physics, 10 (2019) 1454-1486. 
https://doi.org/10.4236/jmp.2019.1012097   \\
\noindent
[40] Pommaret, J.-F.: Homological Solution of the Lanczos Problems in Arbitrary Dimension, Journal of Modern Physics, 12 (2021) 829-858.   https://arxiv.org/abs/1803.09610. \\
\noindent
https://doi.org/10.4236/jmp.2020.1110104  \\
\noindent
[41]  Pommaret, J.-F.: The Mathematical Foundations of Elasticity and Electromagnetism Revisited, Journal of Modern Physics, 10 (2019) 1566-1595.     \\
 https://doi.org/10.4236/jmp.2019.1013104 (https://arxiv.org/abs/1802.02430 ) \\
\noindent
[42] Pommaret, J.-F.: Minimum Parametrization of the Cauchy Stress Operator, Journal of modern Physics, 12 (2021) 453-482. (https://arxiv.org/abs/2101.03959) \\
\noindent 
https://doi.org/10.4236/jmp.2021.124032  \\
\noindent
[43] Pommaret, J.-F.: The Conformal Group Revisited, https://arxiv.org/abs/2006.03449   \\
\noindent
https://doi.org/10.4236/jmp.2021.1213106.  \\
\noindent
[44]  Pommaret, J.-F.: Nonlinear Conformal Electromagnetism, https://arxiv.org/abs/2007.01710 \\
\noindent
https://doi.org/10.4236/jmp.2022.134031   \\
\noindent
[45] Pommaret, J.-F.: Killing Operator for the Kerr Metric, (https://arxiv.org/abs/2211.00064)  \\
\noindent
[46] Pommaret, J.-F.: How Many Structure Constants do Exist in Riemannian Geometry ?, \\
Mathematics in Computer Science, 16, 23 (2022). https://doi.org/10.1007/s11786-022-00546-3   \\
\noindent
[47] Pommaret, J.-F.: Gravitational Waves and Parametrizations of Linear Differential Operators, \\ https://www.intechopen.com/online-first/1119249 
(https://doi.org/10.5992/intechopen.1000851)  \\
\noindent
[48] Pommaret, J.-F.: Gravitational Waves and Lanczos Potentials, JMP, to appear.  \\
\noindent
[49] Pommaret, J.-F.: General Relativity and Gauge Theory: Beyond the Mirror, arXiv: 2302.06585.       \\
\noindent
[50] Pommaret, J.-F., Quadrat, A.: Localization and Parametrization of Linear Multidimensional Control Systems, Systems \& Control Letters, 37 (1999) 247-260.  \\
\noindent 
[51] Robertz, D.: An Involutive GVW Algorithm and the Computation of Pommaret Bases, Mathematics in Computer Science, https://doi.org/10.1007/s11786-021-00512-5  \\
\noindent
[52] Rotman, J.J.: An Introduction to Homological Algebra, Pure and Applied Mathematics, Academic Press (1979).  \\
\noindent
[53] Schneiders, J.-P.: An Introduction to D-Modules, Bull. Soc. Roy. Sci. Li\`{e}ge, 63, 223-295 (1994).  \\
\noindent
[54] Spencer, D.C.: Overdetermined Systems of Partial Differential Equations, Bull. Am. Math. Soc., 75 (1965) 1-114.\\
\noindent
[55] Weyl, H.: Space, Time, Matter, (1918) (Dover, 1952).  \\
\noindent
[56] Zharkov, A.Yu., Blinkov, Yu.A: Involution Approach to investigating Polynomial systems, Mathematics and Computer Simulation, 42, 4 (1996) 323-332. \\
https://doi.org/10.1016/S0378-4754(96)00006-7  \\

\end{document}